\documentclass{amsart}

\usepackage[a4paper,headheight=0.3cm]{geometry}
\usepackage{fancyhdr}
\usepackage{url}
\usepackage{amsmath}
\usepackage{amsthm}
\usepackage{amsfonts}
\usepackage{mathtools}
\usepackage{tikz}
\usepackage{comment}

\usepackage{subfig} 
\usepackage{mathdots} 
\usepackage{enumitem}
\usepackage{url}
\usepackage[pdftitle={Identifying Subtree Perfectness in Decision Trees},%
pdfauthor={Nathan Huntley, Matthias C. M. Troffaes},%
pdfkeywords={choice function, subtree perfectness, backward induction, sequential, decision tree, consequentialist, maximality, maximin, E-admissibility, imprecise probability, imprecise utility}]{hyperref}

\title{Identifying Subtree Perfectness in Decision Trees}
\author{Nathan Huntley}
\address{SYSTeMS Research Group, Universiteit Gent, Belgium}
\email{nath.huntley@gmail.com}
\author{Matthias C. M. Troffaes}
\address{Dept. of Mathematical Sciences, Durham University, UK}
\email{matthias.troffaes@durham.ac.uk}

\makeatletter

\newcommand{\treenotation}[1]{\treenotationevenbetter{#1}}

\newcommand{\treenotationevenbetter}[1]{
\treenotationstart{#1}}

\def\treenotationstart#1{
\@ifnextchar[{
\treenotationsqsq{#1}}%
{
\@ifnextchar({
\treenotationsqro{#1}
}
{
#1}}}%

\def\treenotationsqsq#1[#2]{%
\@ifnextchar[{
\treenotationsqsq{#1{}_{#2}}}%
{
\@ifnextchar({
\treenotationsqro{#1{}_{#2}}
}
{
#1{}_{#2}}}}%

\def\treenotationsqro#1(#2){%
\@ifnextchar[{
\treenotationrosq{#1{}_{#2}}}%
{
\@ifnextchar({
\treenotationroro{#1{}_{#2}}
}
{
#1{}_{#2}}}}%

\def\treenotationroro#1(#2){%
\@ifnextchar[{
\treenotationrosq{#1{}_{#2}}}%
{
\@ifnextchar({
\treenotationroro{#1{}_{#2}}
}
{
#1{}_{#2}}}}%

\def\treenotationrosq#1[#2]{%
\@ifnextchar[{
\treenotationsqsq{#1{}_{#2}}}%
{
\@ifnextchar({
\treenotationsqro{#1{}_{#2}}
}
{
#1{}_{#2}}}}%

\makeatother

\newcommand{\pr}{P}

\newcommand{\lnex}{\underline{E}}
\newcommand{\unex}{\overline{E}}

\newcommand{\pspace}{\Omega}

\newcommand{\domlinprevs}{\mathcal{M}}

\newcommand{\SetR}{\mathbb{R}}

\newcommand{\decnode}{\treenotation{N}}
\newcommand{\chancenode}{\treenotation{N}}

\newcommand{\rewardset}{\mathcal{R}}

\newtheorem{lemma}{Lemma}

\newtheorem{theorem}[lemma]{Theorem}
\newtheorem{example}[lemma]{Example}
\newtheorem{definition}[lemma]{Definition}

\newcommand{\tree}{T}
\newcommand{\atree}{T'}

\newcommand{\compl}[1]{{#1}^c}

\begin{document}

\begin{abstract}
  In decision problems, often, utilities and probabilities are hard to determine.
  In such cases,
  one can resort to so-called choice functions.
  They provide a means to determine which options in a particular set
  are optimal,
  and allow incomparability among any number of options.
  Applying choice functions in sequential decision problems can be highly
  non-trivial, as the usual properties of maximizing expected utility
  may no longer be satisfied.
  In this paper, we study one of these properties:
  we revisit and reinterpret Selten's concept of subgame perfectness
  in the context of decision trees,
  leading us to the concept of subtree perfectness,
  which basically says that the optimal solution of a decision tree
  should not depend on any larger tree it may be embedded in.
  In other words, subtree perfectness excludes counterfactual reasoning,
  and therefore may be desirable
  from some philosophical points of view.
  Subtree perfectness is also desirable from a practical point of view,
  because it admits efficient algorithms
  for solving decision trees, such as backward induction.
  The main contribution of this paper is a very simple non-technical criterion
  for determining whether any given choice function
  will satisfy subtree perfectness or not.
  We demonstrate the theorem
  and illustrate subtree perfectness, or the lack thereof,
  through numerous examples,
  for a wide variety of choice functions,
  where incomparability
  among strategies can be caused by imprecision
  in either probabilities or utilities.
  We find that almost no
  choice function, except for maximizing expected utility, satisfies it
  in general.
  We also find that choice functions other than maximizing expected utility
  can satisfy it, provided that we restrict either the structure of the tree,
  or the structure of the choice function.
\end{abstract}

\maketitle

\thispagestyle{fancy}

\section{Introduction}
\label{sec:intro}

In statistical decision problems,
one is often faced with too little data for too many parameters,
and elicitation of all probabilities and utilities
can be a prohibitively expensive process, rendering full modeling infeasible.
In such cases, it may still be possible to bound the set of reasonable
models, i.e. to bound probability distributions and utility functions
\cite{1854:boole,1931:ramsey,1975:williams:condprev,1980:levi,1983:kyburg,1988:seidenfeld}.
The way in which such models can be used for decision making is well studied
\cite{1991:walley,1999:jaffray,2001:augustin,2002:augustin,2004:seidenfeld,2005:kikuti,2006:jaffray,2007:seidenfeld,2007:troffaes,2011:kikuti:sequential}.

In this paper, we consider a very large class of decision models, namely any that can be represented by conditional choice functions on gambles: we merely assume that for any finite set of gambles (functions from the possibility space $\pspace$ to a set of rewards $\rewardset$; these generalize random variables, or horse lotteries), and any conditioning event, the subject can give a non-empty subset of gambles that he considers optimal if the conditioning event were to occur. Gambles that are non-optimal would never be selected, and the subject is unable to express further preference between the optimal ones. Maximizing conditional expected utility is a simple example of such a choice function.
General choice functions, however, need neither probability nor utility---not even a total preorder, and in fact, not even a partial one, although obviously these obtain as special cases.

Then,
in a single agent sequential decision problem, at any stage, one has two ways of looking at its solution: the problem can be considered either in its simplest form---discarding any past stages, or as part of a much larger problem---possibly considering choices and events that did not actually obtain. A reasonable requirement is that, at any stage, the solution is independent of the larger problem it is embedded in. In this paper, we call this requirement \emph{subtree perfectness}.

Selten~\cite{1975:selten} introduced a similar idea for multi-agent extensive form games, called \emph{subgame perfectness}.
Solutions of extensive form games take the form of \emph{equilibrium points}.
If, for every subgame (a part of a game that is again a game), the restriction of the equilibrium point of the full game to that subgame also yields an equilibrium point of that subgame, then the equilibrium point of the full game is called \emph{subgame perfect}.
Selten showed that for games with perfect recall, (perfect\footnote{Without going into much detail, a \emph{perfect equilibrium point} is one which is stable under small perturbations \cite[p.~38]{1975:selten}.}) equilibrium points are subgame perfect \cite[p.~39, Thm.~2]{1975:selten}, that is, they are independent of any larger game in which they could be embedded.

In this paper we investigate single agent sequential decision making modeled by decision trees. Although such problems differ in many ways from extensive form games, subtree perfectness is clearly analogous to subgame perfectness, as the following example shows.

\begin{figure}
  \hfill
  \begin{minipage}[t]{.45\textwidth}
  \begin{center}
    \begin{tikzpicture}
      [minimum size=2em,parent anchor=east,child anchor=west,grow'=east,transform shape]
      \node[draw,rectangle]{$N_1$}
      [sibling distance=3em, level distance=4em]
      child{
        node[draw,rectangle]{$N_2$}
        [sibling distance=2em]
        child{
          node[right]{cake}
        }
        child{
          node[right]{ice cream}
        }
      }
      child{
        node[right]{scones}
      };
    \end{tikzpicture}
    \caption{Two-stage problem.}
    \label{fig:two:stage:problem}
  \end{center}
  \end{minipage}
  \hfill
  \begin{minipage}[t]{.45\textwidth}
  \begin{center}
    \begin{tikzpicture}[minimum size=2em,parent anchor=east,child anchor=west,grow'=east,transform shape]
      \node[draw,rectangle]{$N_2$}
      [sibling distance=2em]
      child{
        node[right]{cake}
      }
      child{
        node[right]{ice cream}
      };
    \end{tikzpicture}
    \caption{Second stage.}
    \label{fig:second:stage}
  \end{center}
  \end{minipage}
  \hfill
\end{figure}

Consider the decision problem in Fig.~\ref{fig:two:stage:problem} (replicated from~\cite{2011:huntley:subtree:perfectness}). In the first stage, the subject chooses between taking scones, or proceeding to the second stage. In the second stage, the subject chooses between cake, or ice cream. Suppose the subject prefers to reject scones and to choose ice cream at the second stage. This strategy induces a substrategy in the subtree for the second stage: choose ice cream over cake.

But, as with multi-agent games, we can instead consider the subtree for the second stage separately, as in Fig.~\ref{fig:second:stage}. If, in this smaller tree, the subject prefers ice cream, then his solution is \emph{subtree perfect}: his solution for the full tree induces a strategy in the subtree, and this strategy coincides with his solution for the subtree. If the subject states a different preference (either no preference, or clear preference for cake), then his solution lacks subtree perfectness. So, subtree perfectness essentially means that the optimal induced strategies in a subtree do not depend on the full tree in which the subtree is embedded.

Our main result, given in Theorem~\ref{thm:subtreeperfectness} further in the paper, is an extremely straightforward necessary and sufficient condition on choice functions for subtree perfectness to hold in any decision tree.
A similar result was given rigorous treatment in~\cite{2011:huntley:subtree:perfectness}, however that treatment was rather technical and opaque,
and not as obvious to apply to practical problems.
The main contribution of this paper is to provide an overview of the result that is accessible and easy to apply, but still general enough to apply to any choice function.
Therefore, there is far more emphasis on rigorous examples, diagrams, and procedures, and far less on detailed definitions and lengthy proofs.
In so doing, we establish a useful philosophical and practical link between the work that has been done on the behavior of specific choice functions, such as~\cite{1999:jaffray,2002:augustin,2004:seidenfeld}, the more general and theoretical works such as~\cite{1988:hammond,1989:machina,1990:mcclennen,2011:huntley:subtree:perfectness}, and proposed algorithms for solving decision trees~\cite{2002:harmanec,2005:kikuti,2005:decooman,2008:huntley:troffaes::impdectrees:smps,2011:kikuti:sequential,2012:huntley::backinduct}.
Readers seeking formal proofs of the results, or those who are interested in proving similar results about sequential decision problems, are directed to~\cite{2011:huntley:subtree:perfectness}.

In this paper we consider only the traditional normal form method of listing all strategies,
finding their corresponding gambles, applying a choice function, and listing all strategies that induce optimal gambles. In doing so, we assume act-state independence---that is, choice functions do not depend on the decision.
This type of solution is also called \emph{resolute}~\cite{1990:mcclennen,2011:kikuti:sequential}.
Several related works have considered the alternative \emph{extensive form solutions}, such as in Hammond~\cite{1988:hammond}, Seidenfeld~\cite{1988:seidenfeld}, and Machina~\cite{1989:machina}.
We do not consider extensive form solutions in this paper, although the results are closely related.
Without going into details, subtree perfectness is similar to Hammond's \emph{consistency} and \emph{consequentialism}, Machina's \emph{separability over mutually exclusive events}, and McClennen's \emph{separability} and \emph{dynamic consistency}~\cite{1990:mcclennen}.
Further details about the links with extensive form subtree perfectness can be found in\cite{2011:huntley:subtree:perfectness}. 
Subtree perfectness also is related to the weaker concept of \emph{backward induction}.
More details on normal form backward induction can be found in~\cite{2011:huntley:subtree:perfectness,2012:huntley::backinduct}.

It should also be noted that the algorithm proposed by Kikuti et al.~\cite{2005:kikuti,2011:kikuti:sequential} follows exactly our concept of normal form backward induction, although they then use it to construct an extensive form (what they call \emph{consequentialist}) solution.
They note that, for all the choice functions they consider, their solution is ``inconsistent'' with our normal form solution, which essentially means that subtree perfectness is violated.
For two these choice functions, maximality and E-admissibility, it is observed that this inconsistency is not so severe: although arcs that are non-optimal in the full tree may appear in local trees, no arc that is non-optimal in a local tree may reappear in the full tree.
In Fig.~\ref{fig:two:stage:problem}, for example, it would not be possible under these criteria for the subject to say they preferred cake to ice cream at $\decnode[2]$ and then at $\decnode[1]$ include ice cream in the optimal set.
This weaker form of subtree perfectness occurs because the two criteria satisfy the conditions of our backward induction theorem~\cite{2012:huntley::backinduct}.
What this means in practice is that for these criteria the work of Kikuti et al. can be used to efficiently construct the consequentialist solution that they are interested in, and also the normal form solution that we are interested in, and although subtree perfectness does not hold, both the solutions are at least somewhat consistent and sensible.

The paper is structured as follows. Section~\ref{sec:introductory:example} provides an introductory example using expected utility to outline what we mean by a solution of a tree, and what we mean by subtree perfectness.
Sections~\ref{sec:subtree:perfectness:definitions} and~\ref{sec:definition:choice:function} give rigorous, yet non-technical, definitions of the basic concepts required to formulate our main result.
That result, the subtree perfectness theorem, is provided in Section~\ref{sec:theorem}.
We then demonstrate the theorem on numerous choice functions and decision trees: Section~\ref{sec:subtreeperfectness:examples} gives further examples of satisfaction of subtree perfectness for certain choice functions and trees, whilst Section~\ref{sec:failures:subtreeperfectness:examples} gives further examples of failures of subtree perfectness.
Section~\ref{sec:conclusion} concludes the paper.

\section{Subtree Perfectness}\label{sec:subtreeperfectness}

\subsection{Examples: Expected Utility}\label{sec:introductory:example}

\begin{figure}
  \begin{center}
    \begin{tikzpicture}
      [minimum size=2em,parent anchor=east,child anchor=west,grow'=east,transform shape]
      \node[draw,rectangle]{$\decnode$}
      [level distance = 6em, sibling distance=8em]
        child{
          node[draw,rectangle]{$\decnode[1]$}
          [sibling distance=4em]
          child{
            node[draw, circle]{$\chancenode[1](1)$}
            [sibling distance=2em]
            child{
	      node[right]{$2$}
	      edge from parent
	      node[above, sloped]{$E_1$}
            }
            child{
	      node[right]{$0$}
	      edge from parent
	      node[below, sloped]{$E_2$}
            }
          }
          child{
            node[draw, circle]{$\chancenode[1](2)$}
            [sibling distance=2em]
            child{
	      node[right]{$1$}
	      edge from parent
	      node[above, sloped]{$E_1$}
            }
            child{
	      node[right]{$1$}
	      edge from parent
	      node[below, sloped]{$E_2$}
            }
          }
        }
        child{
          node[draw,rectangle]{$\decnode[2]$}
          [sibling distance=4em]
          child{
            node[draw, circle]{$\chancenode[2](1)$}
            [sibling distance=2em]
            child{
	      node[right]{$3$}
	      edge from parent
	      node[above, sloped]{$E_1$}
            }
            child{
	      node[right]{$0$}
	      edge from parent
	      node[below, sloped]{$E_2$}
            }
          }
          child{
            node[draw, circle]{$\chancenode[2](2)$}
            [sibling distance=2em]
            child{
	      node[right]{$4$}
	      edge from parent
	      node[above, sloped]{$E_1$}
            }
            child{
	      node[right]{$-1$}
	      edge from parent
	      node[below, sloped]{$E_2$}
            }
          }
        };        
    \end{tikzpicture}
  \end{center}
  \caption{Decision tree for first expected utility example.}
  \label{fig:expected:utility:one}
\end{figure}
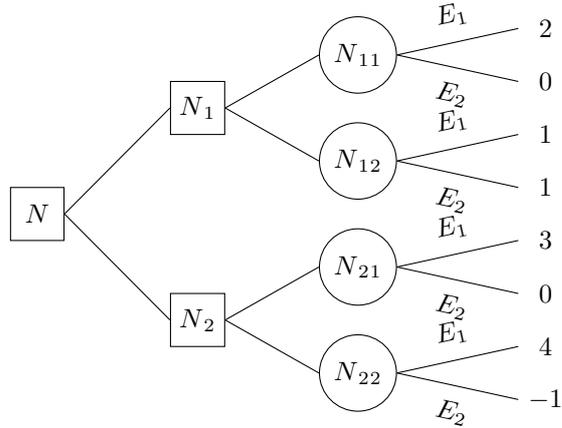

We begin by considering some examples of solving small decision trees using expected utility.
Of course, these solutions are almost trivial to find, but we spend some time on them to outline our ideas in a familiar environment.
Consider the decision tree in Fig.~\ref{fig:expected:utility:one}, and suppose that $\pr(E_1)=0.6$, $\pr(E_2)=0.4$.
Square \emph{decision} nodes are points where the subject chooses the arc to follow.
Round \emph{chance} nodes are points where the arc followed depends on the (initially unknown) state of nature, in this case $E_1$ or $E_2$.
Typically one would approach the problem by \emph{backward induction} \cite{1985:lindley}, but we consider an alternative method.
This is because for more general choice functions, backward induction is more complicated and may not give the desired answer.

Consider a path that the subject may take through the tree.
At $N$, she has a choice of arcs and suppose she follows the arc to $\decnode[1]$.
Here again she has a choice; suppose she takes the arc to $\decnode[1](1)$.
Here, she has no choice over the arc followed: nature decides.
Thus, for this combination of her choices, her final reward is $2$ if $E_1$ occurs and $0$ otherwise.
We can calculate the expected utility of this \emph{strategy} to be $1.2$.

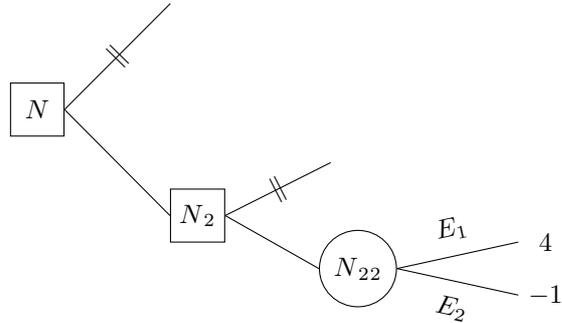
\begin{figure}
  \begin{center}
    \begin{tikzpicture}
      [minimum size=2em,parent anchor=east,child anchor=west,grow'=east,transform shape]
      \node[draw,rectangle]{$\decnode$}
      [level distance = 6em, sibling distance=8em]
        child{
          node{}
          edge from parent
          node[sloped]{$\|$}
        }
        child{
          node[draw,rectangle]{$\decnode[2]$}
          [sibling distance=4em]
          child{
            node{}
            edge from parent
            node[sloped]{$\|$}
          }
          child{
            node[draw, circle]{$\chancenode[2](2)$}
            [sibling distance=2em]
            child{
	      node[right]{$4$}
	      edge from parent
	      node[above, sloped]{$E_1$}
            }
            child{
	      node[right]{$-1$}
	      edge from parent
	      node[below, sloped]{$E_2$}
            }
          }
        };        
    \end{tikzpicture}
  \end{center}
  \caption{Normal form solution for first expected utility example.}
  \label{fig:expected:utility:one:normal:form:solution}
\end{figure}

There are three other strategies in the tree: the paths $\decnode\to\decnode[1]\to\chancenode[1](2)$, $\decnode\to\decnode[2]\to\chancenode[2](1)$, and $\decnode\to\decnode[2]\to\chancenode[2](2)$.
These have expected utilities $1$, $1.8$, and $2$ respectively.
We therefore conclude that the solution of the tree is to take the path from $\decnode$ to $\decnode[2]$ and then to $\chancenode[2](2)$.
This is illustrated in the decision tree in Fig.~\ref{fig:expected:utility:one:normal:form:solution}.

\begin{figure}
    \begin{center}
      \begin{tikzpicture}
        [minimum size=2em,parent anchor=east,child anchor=west,grow'=east,transform shape]
        \node[draw,rectangle]{$\decnode[1]$}
          [sibling distance=4em]
          child{
            node[draw, circle]{$\chancenode[1](1)$}
            [sibling distance=2em]
            child{
	      node[right]{$2$}
	      edge from parent
	      node[above, sloped]{$E_1$}
            }
            child{
	      node[right]{$0$}
	      edge from parent
	      node[below, sloped]{$E_2$}
            }
          }
          child{
            node[draw, circle]{$\chancenode[1](2)$}
            [sibling distance=2em]
            child{
	      node[right]{$1$}
	      edge from parent
	      node[above, sloped]{$E_1$}
            }
            child{
	      node[right]{$1$}
	      edge from parent
	      node[below, sloped]{$E_2$}
            }
          };        
      \end{tikzpicture}
      \hspace{2em}
    \begin{tikzpicture}
      [minimum size=2em,parent anchor=east,child anchor=west,grow'=east,transform shape]
      \node[draw,rectangle]{$\decnode[2]$}
          [sibling distance=4em]
          child{
            node[draw, circle]{$\chancenode[2](1)$}
            [sibling distance=2em]
            child{
	      node[right]{$3$}
	      edge from parent
	      node[above, sloped]{$E_1$}
            }
            child{
	      node[right]{$0$}
	      edge from parent
	      node[below, sloped]{$E_2$}
            }
          }
          child{
            node[draw, circle]{$\chancenode[2](2)$}
            [sibling distance=2em]
            child{
	      node[right]{$4$}
	      edge from parent
	      node[above, sloped]{$E_1$}
            }
            child{
	      node[right]{$-1$}
	      edge from parent
	      node[below, sloped]{$E_2$}
            }
          };        
    \end{tikzpicture}
    \caption{Subtrees for first expected utility example.}
    \label{fig:expected:utility:one:subtrees}
  \end{center}
\end{figure}
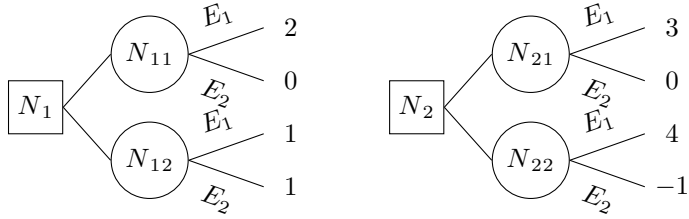

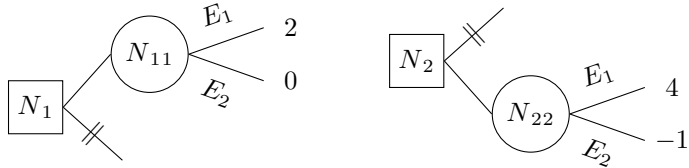
\begin{figure}
    \begin{center}
      \begin{tikzpicture}
        [minimum size=2em,parent anchor=east,child anchor=west,grow'=east,transform shape]
        \node[draw,rectangle]{$\decnode[1]$}
          [sibling distance=4em]
          child{
            node[draw, circle]{$\chancenode[1](1)$}
            [sibling distance=2em]
            child{
	      node[right]{$2$}
	      edge from parent
	      node[above, sloped]{$E_1$}
            }
            child{
	      node[right]{$0$}
	      edge from parent
	      node[below, sloped]{$E_2$}
            }
          }
          child{
            node{}
	      edge from parent
	      node[sloped]{$\|$}
            };        
      \end{tikzpicture}
      \hspace{2em}
    \begin{tikzpicture}
      [minimum size=2em,parent anchor=east,child anchor=west,grow'=east,transform shape]
      \node[draw,rectangle]{$\decnode[2]$}
          [sibling distance=4em]
          child{
            node{}
	      edge from parent
	      node[sloped]{$\|$}
            }
          child{
            node[draw, circle]{$\chancenode[2](2)$}
            [sibling distance=2em]
            child{
	      node[right]{$4$}
	      edge from parent
	      node[above, sloped]{$E_1$}
            }
            child{
	      node[right]{$-1$}
	      edge from parent
	      node[below, sloped]{$E_2$}
            }
          };        
    \end{tikzpicture}
    \caption{Solutions of subtrees for first expected utility example.}
    \label{fig:expected:utility:one:subtrees:solution}
  \end{center}
\end{figure}

Now, subtree perfectness is a property relating solutions of large trees to solutions of subtrees.
In this particular example, there are two non-trivial subtrees of the full tree, shown in Fig.~\ref{fig:expected:utility:one:subtrees}.
We can solve these trees in the same manner; indeed we have already calculated the required expected utilities.
The optimal paths are $\decnode[1]\to\chancenode[1](1)$ in the left-hand tree,
and $\decnode[2]\to\chancenode[2](2)$ in the right-hand tree.
These solutions are displayed in Fig~\ref{fig:expected:utility:one:subtrees:solution}.

Are these local solutions consistent with the solution in Fig.~\ref{fig:expected:utility:one:normal:form:solution}?
The first thing to note is that $\decnode[1]$ does not appear at all in Fig.~\ref{fig:expected:utility:one:normal:form:solution}.
We can therefore ignore it: behavior in a subtree that will never actually be reached is irrelevant to subtree perfectness.
At $\decnode[2]$, we see that both the solutions (Fig.~\ref{fig:expected:utility:one:normal:form:solution} and Fig.~\ref{fig:expected:utility:one:subtrees:solution}) contain only the path to $\chancenode[2](2)$, so there is no inconsistency here either.
We say that subtree perfectness holds for this problem.

\begin{figure}
  \begin{center}
    \begin{tikzpicture}
      [minimum size=2em,parent anchor=east,child anchor=west,grow'=east,transform shape]
      \node[draw,circle]{$\chancenode$}
      [level distance = 6em, sibling distance=8em]
        child{
          node[draw,rectangle]{$\decnode[1]$}
          [sibling distance=4em]
          child{
            node[draw, circle]{$\chancenode[1](1)$}
            [sibling distance=2em]
            child{
	      node[right]{$2$}
	      edge from parent
	      node[above, sloped]{$E_1$}
            }
            child{
	      node[right]{$0$}
	      edge from parent
	      node[below, sloped]{$E_2$}
            }
          }
          child{
            node[draw, circle]{$\chancenode[1](2)$}
            [sibling distance=2em]
            child{
	      node[right]{$1$}
	      edge from parent
	      node[above, sloped]{$E_1$}
            }
            child{
	      node[right]{$1$}
	      edge from parent
	      node[below, sloped]{$E_2$}
            }
          }
          edge from parent
          node[above,sloped]{$A_1$}
        }
        child{
            node[draw, circle]{$\chancenode(2)$}
            [sibling distance=4em]
            child{
	      node[right]{$3$}
	      edge from parent
	      node[above, sloped]{$E_1$}
            }
            child{
	      node[right]{$0$}
	      edge from parent
	      node[below, sloped]{$E_2$}
            }
            edge from parent
            node[below,sloped]{$A_2$}
          };        
    \end{tikzpicture}
  \end{center}
  \caption{Decision tree for second and third expected utility examples.}
  \label{fig:expected:utility:two}
\end{figure}
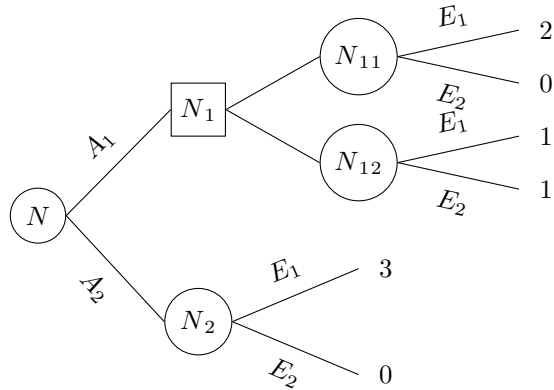

Next we move onto a slightly more involved example, in which the initial node is a chance node, as shown in Fig.~\ref{fig:expected:utility:two}.
We need a few more probabilities for this example: let $\pr(E_1|A_1)=0.6$, $\pr(E_2|A_1)=0.4$, $\pr(E_1|A_2)=0.4$, $\pr(E_2|A_2)=0.6$, $\pr(A_1)=0.5$, and $\pr(A_2)=0.5$.

\begin{figure}
  \begin{center}
    \begin{tikzpicture}
      [minimum size=2em,parent anchor=east,child anchor=west,grow'=east,transform shape]
      \node[draw,circle]{$\chancenode$}
      [level distance = 6em, sibling distance=8em]
        child{
          node[draw,rectangle]{$\decnode[1]$}
          [sibling distance=4em]
          child{
            node[draw, circle]{$\chancenode[1](1)$}
            [sibling distance=2em]
            child{
	      node[right]{$2$}
	      edge from parent
	      node[above, sloped]{$E_1$}
            }
            child{
	      node[right]{$0$}
	      edge from parent
	      node[below, sloped]{$E_2$}
            }
          }
          child{
            node{}
            edge from parent
            node[sloped]{$\|$}
          }
          edge from parent
          node[above,sloped]{$A_1$}
        }
        child{
            node[draw, circle]{$\chancenode(2)$}
            [sibling distance=4em]
            child{
	      node[right]{$3$}
	      edge from parent
	      node[above, sloped]{$E_1$}
            }
            child{
	      node[right]{$0$}
	      edge from parent
	      node[below, sloped]{$E_2$}
            }
            edge from parent
            node[below,sloped]{$A_2$}
          };        
    \end{tikzpicture}
  \end{center}
  \caption{Normal form solution for second expected utility example.}
  \label{fig:expected:utility:two:normal:form:solution}
\end{figure}
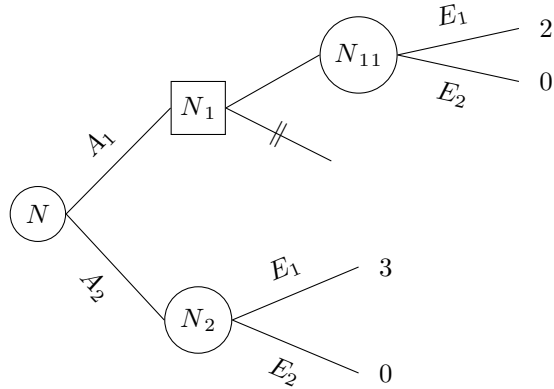

In this example, the two available strategies are ``if $A_1$, take the path $\decnode[1]\to\chancenode[1](1)$'', and ``if $A_1$, take the path $\decnode[1]\to\chancenode[1](2)$''.
Note that there is no decision node following $A_2$, so the subject does not need to specify anything in that case.
These strategies have slightly more complicated gambles than before.
Consider the first: if $A_1$ and $E_1$ both occur, the outcome is $2$; if $A_1$ and $E_2$ both occur, the outcome is $0$; if $A_2$ and $E_1$ both occur, the outcome is 3; and if $A_2$ and $E_2$ both occur, the outcome is $0$.
This gives expected utility $0.5 (0.6 \cdot 2 + 0.4 \cdot 0) + 0.5 ( 0.4 \cdot 3 + 0.6 \cdot 0)=1.2$.
The second strategy gives expected utility $0.5\cdot 1 + 0.5 ( 0.4 \cdot 3 + 0.6 \cdot 0)=1.1$.
Hence the solution is the first strategy, shown in Fig~\ref{fig:expected:utility:two:normal:form:solution}.

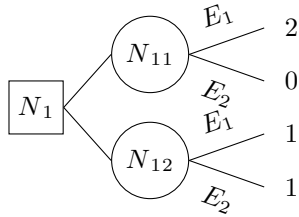
\begin{figure}
  \begin{center}
    \begin{tikzpicture}
      [minimum size=2em,parent anchor=east,child anchor=west,grow'=east,transform shape]
          \node[draw,rectangle]{$\decnode[1]$}
          [sibling distance=4em]
          child{
            node[draw, circle]{$\chancenode[1](1)$}
            [sibling distance=2em]
            child{
	      node[right]{$2$}
	      edge from parent
	      node[above, sloped]{$E_1$}
            }
            child{
	      node[right]{$0$}
	      edge from parent
	      node[below, sloped]{$E_2$}
            }
          }
          child{
            node[draw, circle]{$\chancenode[1](2)$}
            [sibling distance=2em]
            child{
	      node[right]{$1$}
	      edge from parent
	      node[above, sloped]{$E_1$}
            }
            child{
	      node[right]{$1$}
	      edge from parent
	      node[below, sloped]{$E_2$}
            }
          };        
    \end{tikzpicture}
  \end{center}
  \caption{Subtree for second expected utility example.}
  \label{fig:expected:utility:two:subtrees}
\end{figure}

Again, we look for non-trivial subtrees.
In this example, there is only one, at $\decnode[1]$.
This is shown in Fig~\ref{fig:expected:utility:two:subtrees}---note that this is exactly the left-hand tree of Fig.~\ref{fig:expected:utility:one:subtrees}.
We can solve this as usual, although here we have to remember that the event $A_1$ has occurred, and so we must use the suitable conditional probabilities.
The expected utilities of the two strategies in this case are $0.6\cdot 2 + 0.4\cdot 0=1.2$ and $1$ respectively, so the normal form solution of this subtree is to take the path $\decnode[1]\to\chancenode[1](1)$ (the same as the left-hand tree in Fig.~\ref{fig:expected:utility:one:subtrees:solution}).

Once again we check for subtree perfectness. In both the global and the local solutions, the only path from $\decnode[1]$ is the one to $\chancenode[1](1)$, hence subtree perfectness is satisfied in this example.

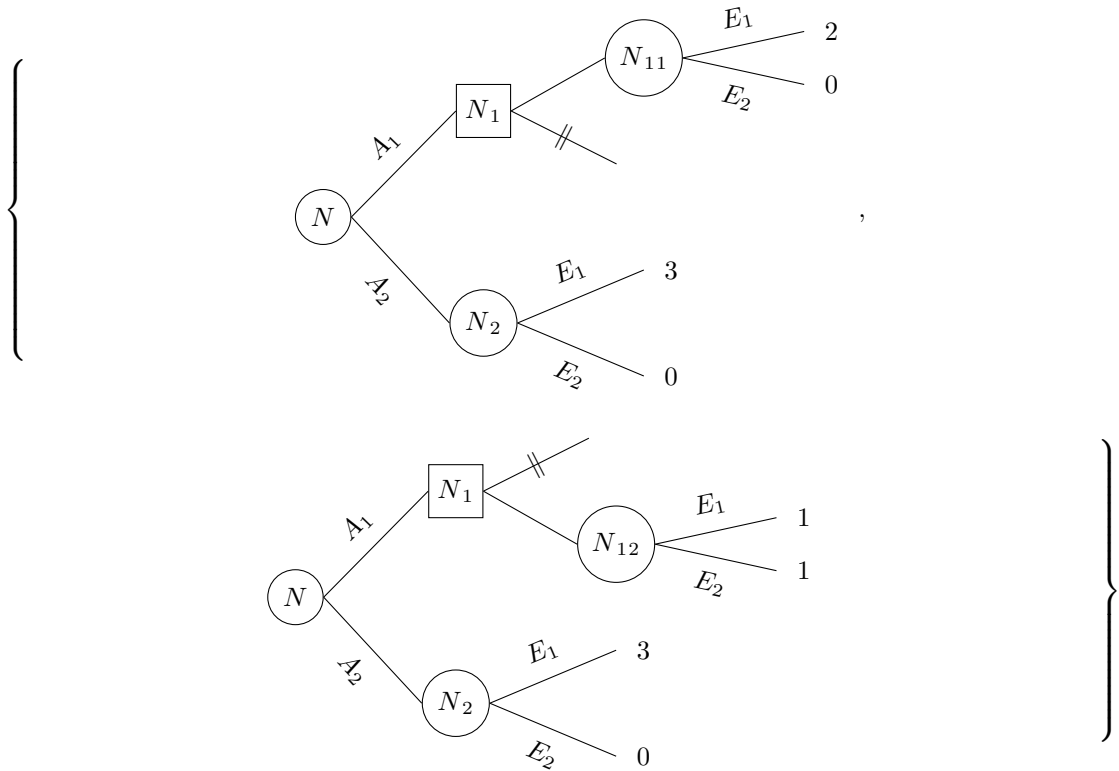
\begin{figure}
  \begin{center}
  $\left\{\rule{0em}{6em}\right.$
    \hfill
    \begin{tikzpicture}
      [minimum size=2em,parent anchor=east,child anchor=west,grow'=east,transform shape,baseline=0em]
      \node[draw,circle]{$\chancenode$}
      [level distance = 6em, sibling distance=8em]
        child{
          node[draw,rectangle]{$\decnode[1]$}
          [sibling distance=4em]
          child{
            node[draw, circle]{$\chancenode[1](1)$}
            [sibling distance=2em]
            child{
	      node[right]{$2$}
	      edge from parent
	      node[above, sloped]{$E_1$}
            }
            child{
	      node[right]{$0$}
	      edge from parent
	      node[below, sloped]{$E_2$}
            }
          }
          child{
            node{}
            edge from parent
            node[sloped]{$\|$}
          }
          edge from parent
          node[above,sloped]{$A_1$}
        }
        child{
            node[draw, circle]{$\chancenode(2)$}
            [sibling distance=4em]
            child{
	      node[right]{$3$}
	      edge from parent
	      node[above, sloped]{$E_1$}
            }
            child{
	      node[right]{$0$}
	      edge from parent
	      node[below, sloped]{$E_2$}
            }
            edge from parent
            node[below,sloped]{$A_2$}
          };        
    \end{tikzpicture},\hfill\mbox{} \\
    \mbox{}\hfill
    \begin{tikzpicture}
      [minimum size=2em,parent anchor=east,child anchor=west,grow'=east,transform shape,baseline=0em]
      \node[draw,circle]{$\chancenode$}
      [level distance = 6em, sibling distance=8em]
        child{
          node[draw,rectangle]{$\decnode[1]$}
          [sibling distance=4em]
          child{
            node{}
            edge from parent
            node[sloped]{$\|$}
          }
          child{
            node[draw, circle]{$\chancenode[1](2)$}
            [sibling distance=2em]
            child{
	      node[right]{$1$}
	      edge from parent
	      node[above, sloped]{$E_1$}
            }
            child{
	      node[right]{$1$}
	      edge from parent
	      node[below, sloped]{$E_2$}
            }
          }
          edge from parent
          node[above,sloped]{$A_1$}
        }
        child{
            node[draw, circle]{$\chancenode(2)$}
            [sibling distance=4em]
            child{
	      node[right]{$3$}
	      edge from parent
	      node[above, sloped]{$E_1$}
            }
            child{
	      node[right]{$0$}
	      edge from parent
	      node[below, sloped]{$E_2$}
            }
            edge from parent
            node[below,sloped]{$A_2$}
          };        
    \end{tikzpicture}\hfill
    $\left.\rule{0em}{6em}\right\}$
  \end{center}
  \caption{Normal form solution for third expected utility example.}
  \label{fig:expected:utility:three:normal:form:solution}
\end{figure}

Finally, we provide an example of expected utility failing subtree perfectness.
As we show later, this requires a conditioning event with probability zero.
Changing the probabilities in the previous example to $\pr(A_1)=0$, $\pr(A_2)=1$ will be sufficient to fail subtree perfectness.
Now the expected utilities of the two strategies in the tree of Fig.~\ref{fig:expected:utility:two} are $0 (0.6 \cdot 2 + 0.4 \cdot 0) + 1 ( 0.4 \cdot 3 + 0.6 \cdot 0)=1.2$ and $0\cdot 1 + 1 ( 0.4 \cdot 3 + 0.6 \cdot 0)=1.2$ respectively.
Hence, we cannot choose between the two strategies using expected utility, and the normal form solution must include both.
This is illustrated in Fig~\ref{fig:expected:utility:three:normal:form:solution}.

Because we did not change any of the conditional probabilities, the normal form solution of the subtree at $\decnode[1]$ will be the same as in the previous example: the left-hand tree of Fig.~\ref{fig:expected:utility:one:subtrees:solution}.
Now we see inconsistency between the two solutions: one of the strategies in the global solution contains the path $\decnode[1]\to\decnode[1](2)$, but this path does not appear in the local solution.
We say that subtree perfectness does not hold at $\decnode[1]$.

Having seen how subtree perfectness works for a standard choice function and simple trees, we now explain how to carry out the steps above for any tree and any choice function.

\subsection{Definition of Subtree Perfectness}\label{sec:subtree:perfectness:definitions}

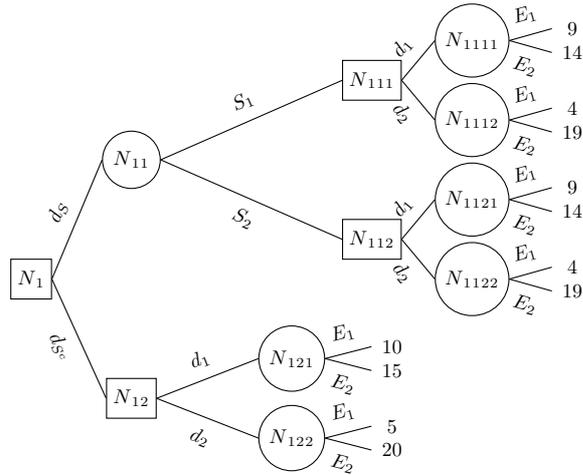
\begin{figure}
  \begin{center}
    \begin{tikzpicture}
      [minimum size=2em,parent anchor=east,child anchor=west,grow'=east,scale=0.75,transform shape]
        \node[draw,rectangle]{$\decnode[1]$}
        [sibling distance=12em,level distance=5em]
        child{
          node[draw,circle]{$\chancenode[1](1)$}
          [sibling distance=8em,level distance=12em]
          child{
            node[draw,rectangle]{$\decnode[1](1)[1]$}
            [sibling distance=4em,level distance=5em]
            child{
              node[draw,circle]{$\chancenode[1](1)[1](1)$}
              [sibling distance=1.2em,level distance=5em]
              child{
                node{$9$}
                edge from parent
                node[above,sloped]{$E_1$}
              }
              child{
                node{$14$}
                edge from parent
                node[below,sloped]{$E_2$}
              }
              edge from parent
              node[above,sloped]{$d_1$}
            }
            child{
              node[draw,circle]{$\chancenode[1](1)[1](2)$}
              [sibling distance=1.2em,level distance=5em]
              child{
                node{$4$}
                edge from parent
                node[above,sloped]{$E_1$}
              }
              child{
                node{$19$}
                edge from parent
                node[below,sloped]{$E_2$}
              }
              edge from parent
              node[below,sloped]{$d_2$}
            }
            edge from parent
            node[above,sloped]{$S_1$}
          }
          child{
            node[draw,rectangle]{$\decnode[1](1)[2]$}
            [sibling distance=4em,level distance=5em]
            child{
              node[draw,circle]{$\chancenode[1](1)[2](1)$}
              [sibling distance=1.2em,level distance=5em]
              child{
                node{$9$}
                edge from parent
                node[above,sloped]{$E_1$}
              }
              child{
                node{$14$}
                edge from parent
                node[below,sloped]{$E_2$}
              }
              edge from parent
              node[above,sloped]{$d_1$}
            }
            child{
              node[draw,circle]{$\chancenode[1](1)[2](2)$}
              [sibling distance=1.2em,level distance=5em]
              child{
                node{$4$}
                edge from parent
                node[above,sloped]{$E_1$}
              }
              child{
                node{$19$}
                edge from parent
                node[below,sloped]{$E_2$}
              }
              edge from parent
              node[below,sloped]{$d_2$}
            }
            edge from parent
            node[below,sloped]{$S_2$}
          }
          edge from parent
          node[above,sloped]{$d_S$}
        }
        child{
          node[draw,rectangle]{$\decnode[1][2]$}
          [sibling distance=4em,level distance=8em]
          child{
            node[draw,circle]{$\chancenode[1][2](1)$}
            [sibling distance=1.2em,level distance=5em]
            child{
              node{$10$}
              edge from parent
              node[above,sloped]{$E_1$}
            }
            child{
              node{$15$}
              edge from parent
              node[below,sloped]{$E_2$}
            }
            edge from parent
            node[above,sloped]{$d_1$}
          }
          child{
            node[draw,circle]{$\chancenode[1][2](2)$}
            [sibling distance=1.2em,level distance=5em]
            child{
              node{$5$}
              edge from parent
              node[above,sloped]{$E_1$}
            }
            child{
              node{$20$}
              edge from parent
              node[below,sloped]{$E_2$}
            }
            edge from parent
            node[below,sloped]{$d_2$}
          }
          edge from parent
          node[below,sloped]{$d_{\compl{S}}$}
        };
    \end{tikzpicture}
  \end{center}
    \caption{A more complex example of a decision tree.}
  \label{fig:lake:tree:basic}
\end{figure}

To state our main result, we need to define subtree perfectness for
arbitrarily large (but finite) decision trees.
In this paper, we will describe and define the concept of subtree perfectness
directly on decision trees, without further formal mathematical notation.
The tree in Fig.~\ref{fig:lake:tree:basic} will serve as a leading example.
Of course, subtree perfectness can also be defined more formally,
using mathematical notation that is very useful for proofs,
as we did in~\cite{2011:huntley:subtree:perfectness},
however such approach is rather technical,
and unnecessary for illustrating the main ideas.
Any reader who is interested
in rigorous mathematical formulations and proofs of properties
relating to choice functions and decision trees
is referred to~\cite{2011:huntley:subtree:perfectness}.
It should also be noted that, although throughout this paper almost all our examples feature real-valued rewards, the language and results we introduce here will work for any form of rewards.

The solutions we consider are sets of \emph{strategies}.
A strategy corresponds to the subject initially specifying all her actions in all eventualities, that is, an action for each decision node she may possibly reach.
Upon arriving at a decision node, she then follows the mandated action.
Thus, after the strategy is chosen, the subject no longer has any control.
In other words, she is a \emph{resolute} decision maker \cite{1990:mcclennen}.

\begin{figure}
  \begin{center}
    \begin{tikzpicture}
      [minimum size=2em,parent anchor=east,child anchor=west,grow'=east,scale=0.75,transform shape]
        \node[draw,rectangle]{$\decnode[1]$}
        [sibling distance=12em,level distance=5em]
        child{
          node[draw,circle]{$\chancenode[1](1)$}
          [sibling distance=8em,level distance=12em]
          child{
            node[draw,rectangle]{$\decnode[1](1)[1]$}
            [sibling distance=4em,level distance=5em]
            child{
              node[draw,circle]{$\chancenode[1](1)[1](1)$}
              [sibling distance=1.2em,level distance=5em]
              child{
                node{$9$}
                edge from parent
                node[above,sloped]{$E_1$}
              }
              child{
                node{$14$}
                edge from parent
                node[below,sloped]{$E_2$}
              }
              edge from parent
              node[above,sloped]{$d_1$}
            }
            child{
              node{}
              edge from parent
              node[below,sloped]{$d_2$}
              node[sloped]{$\|$}
            }
            edge from parent
            node[above,sloped]{$S_1$}
          }
          child{
            node[draw,rectangle]{$\decnode[1](1)[2]$}
            [sibling distance=4em,level distance=5em]
            child{
              node{}
              edge from parent
              node[above,sloped]{$d_1$}
              node[sloped]{$\|$}
            }
            child{
              node[draw,circle]{$\chancenode[1](1)[2](2)$}
              [sibling distance=1.2em,level distance=5em]
              child{
                node{$4$}
                edge from parent
                node[above,sloped]{$E_1$}
              }
              child{
                node{$19$}
                edge from parent
                node[below,sloped]{$E_2$}
              }
              edge from parent
              node[below,sloped]{$d_2$}
            }
            edge from parent
            node[below,sloped]{$S_2$}
          }
          edge from parent
          node[above,sloped]{$d_S$}
        }
        child{
          node{}
          edge from parent
          node[below,sloped]{$d_{\compl{S}}$}
          node[sloped]{$\|$}
        };
    \end{tikzpicture}
  \end{center}
    \caption{One of the strategies of the tree in Fig.~\ref{fig:lake:tree:basic}.}
  \label{fig:strategy:example:tree}
\end{figure}
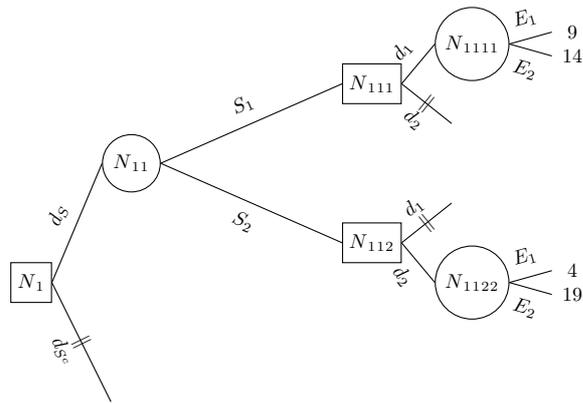

A strategy is easy to represent as a subtree of the initial tree.
Consider the decision tree in Fig.~\ref{fig:lake:tree:basic}, and the strategy ``at $\decnode[1]$ choose $d_S$, then if $S_1$ occurs choose $d_1$ and if $S_2$ occurs choose $d_2$''. This can be represented by the decision tree in Fig.~\ref{fig:strategy:example:tree}.
In this figure, an arc with a double line through means the arc has been deleted.
So the decision tree representation of a strategy is the original tree minus all decision arcs that are not in the strategy.
Note that the labels of the nodes were retained in this process: this will be an important point later on.
Also note that this is now a trivial decision tree, since there is only one option available at each decision node.

In general, we can construct a strategy for a given decision tree by the following method:
\begin{enumerate}
\item For each decision node, pick exactly one decision arc and delete all others.
\item Delete everything not connected to the original root node of the tree.
\end{enumerate}
By applying this method for all possible combinations of choices of decision arcs, we arrive at the set of all available strategies for a given decision tree.\footnote{%
Of course, in practice, there is a more efficient way to
find all strategies---for instance,
in Fig.~\ref{fig:lake:tree:basic}
there are $2^4=16$ combinations of arc choices,
but these lead to only $6$ distinct strategies:
evidently, one does not need to iterate over arcs
of decision nodes in subtrees following deleted arcs.}

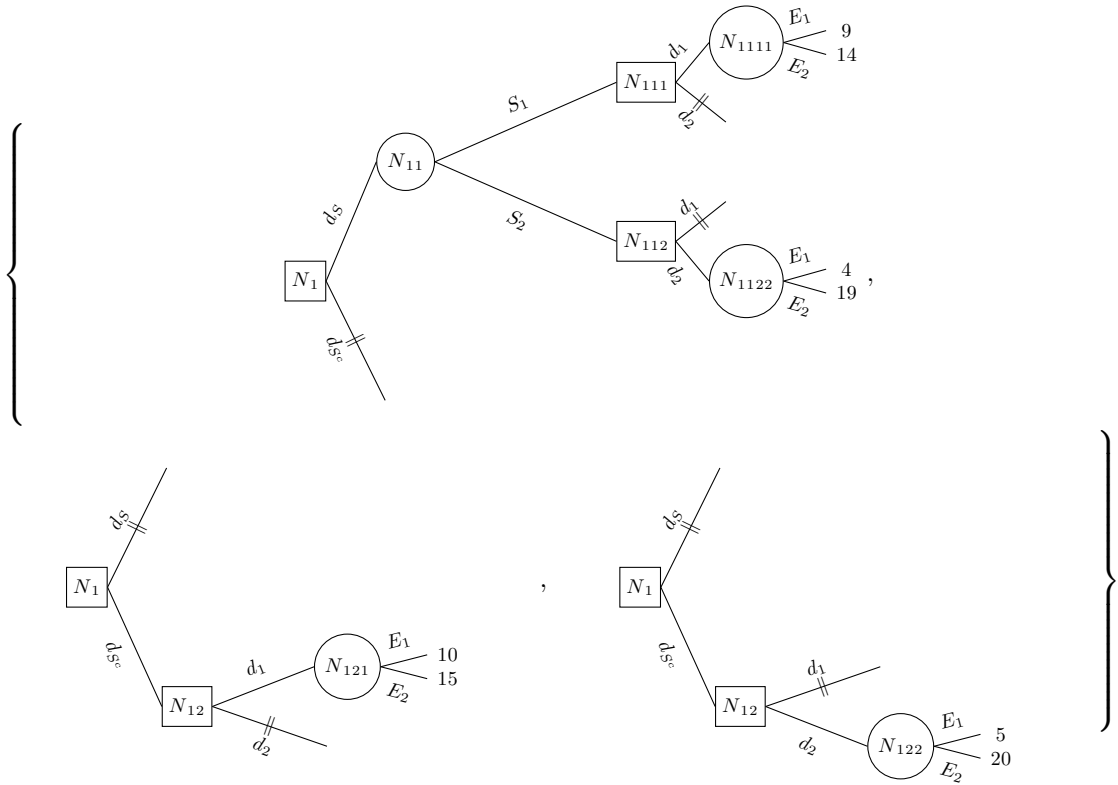
\begin{figure}
  \begin{center}
    $\left\{\rule{0em}{6em}\right.$
    \hfill
    \begin{tikzpicture}
      [minimum size=2em,parent anchor=east,child anchor=west,grow'=east,scale=0.75,transform shape,baseline=0em]
        \node[draw,rectangle]{$\decnode[1]$}
        [sibling distance=12em,level distance=5em]
        child{
          node[draw,circle]{$\chancenode[1](1)$}
          [sibling distance=8em,level distance=12em]
          child{
            node[draw,rectangle]{$\decnode[1](1)[1]$}
            [sibling distance=4em,level distance=5em]
            child{
              node[draw,circle]{$\chancenode[1](1)[1](1)$}
              [sibling distance=1.2em,level distance=5em]
              child{
                node{$9$}
                edge from parent
                node[above,sloped]{$E_1$}
              }
              child{
                node{$14$}
                edge from parent
                node[below,sloped]{$E_2$}
              }
              edge from parent
              node[above,sloped]{$d_1$}
            }
            child{
              node{}
              edge from parent
              node[below,sloped]{$d_2$}
              node[sloped]{$\|$}
            }
            edge from parent
            node[above,sloped]{$S_1$}
          }
          child{
            node[draw,rectangle]{$\decnode[1](1)[2]$}
            [sibling distance=4em,level distance=5em]
            child{
              node{}
              edge from parent
              node[above,sloped]{$d_1$}
              node[sloped]{$\|$}
            }
            child{
              node[draw,circle]{$\chancenode[1](1)[2](2)$}
              [sibling distance=1.2em,level distance=5em]
              child{
                node{$4$}
                edge from parent
                node[above,sloped]{$E_1$}
              }
              child{
                node{$19$}
                edge from parent
                node[below,sloped]{$E_2$}
              }
              edge from parent
              node[below,sloped]{$d_2$}
            }
            edge from parent
            node[below,sloped]{$S_2$}
          }
          edge from parent
          node[above,sloped]{$d_S$}
        }
        child{
          node{}
          edge from parent
          node[below,sloped]{$d_{\compl{S}}$}
          node[sloped]{$\|$}
        };
    \end{tikzpicture},\hfill\mbox{} \\
    \mbox{}\hfill
    \begin{tikzpicture}
      [minimum size=2em,parent anchor=east,child anchor=west,grow'=east,scale=0.75,transform shape,baseline=0em]
        \node[draw,rectangle]{$\decnode[1]$}
        [sibling distance=12em,level distance=5em]
        child{
          node{}
          edge from parent
          node[above,sloped]{$d_S$}
          node[sloped]{$\|$}
        }
        child{
          node[draw,rectangle]{$\decnode[1][2]$}
          [sibling distance=4em,level distance=8em]
          child{
            node[draw,circle]{$\chancenode[1][2](1)$}
            [sibling distance=1.2em,level distance=5em]
            child{
              node{$10$}
              edge from parent
              node[above,sloped]{$E_1$}
            }
            child{
              node{$15$}
              edge from parent
              node[below,sloped]{$E_2$}
            }
            edge from parent
            node[above,sloped]{$d_1$}
          }
          child{
            node{}
            edge from parent
            node[below,sloped]{$d_2$}
            node[sloped]{$\|$}
          }
          edge from parent
          node[below,sloped]{$d_{\compl{S}}$}
        };
    \end{tikzpicture}
    \hfill, \hfill
    \begin{tikzpicture}
      [minimum size=2em,parent anchor=east,child anchor=west,grow'=east,scale=0.75,transform shape,baseline=0em]
        \node[draw,rectangle]{$\decnode[1]$}
        [sibling distance=12em,level distance=5em]
        child{
          node{}
          edge from parent
          node[above,sloped]{$d_S$}
          node[sloped]{$\|$}
        }
        child{
          node[draw,rectangle]{$\decnode[1][2]$}
          [sibling distance=4em,level distance=8em]
          child{
            node{}
            edge from parent
            node[above,sloped]{$d_1$}
            node[sloped]{$\|$}
          }
          child{
            node[draw,circle]{$\chancenode[1][2](2)$}
            [sibling distance=1.2em,level distance=5em]
            child{
              node{$5$}
              edge from parent
              node[above,sloped]{$E_1$}
            }
            child{
              node{$20$}
              edge from parent
              node[below,sloped]{$E_2$}
            }
            edge from parent
            node[below,sloped]{$d_2$}
          }
          edge from parent
          node[below,sloped]{$d_{\compl{S}}$}
        };
    \end{tikzpicture}
    \hfill
    $\left.\rule{0em}{6em}\right\}$
  \end{center}
    \caption{Possible normal form solution of Fig.~\ref{fig:lake:tree:basic}.}
  \label{fig:normal:form:solution:example}
\end{figure}

As in the examples earlier, we consider solutions of decision trees of the following type: first we find all available strategies, then we apply some criterion to them that returns a non-empty subset of strategies.
Ideally, we might want our criterion to return a single strategy.
However, as argued in the introduction, in many situations the subject
may not have enough information to make such a judgment,
for instance because their probabilities and utilities
are not fully determined.
A function that is applied to a set of options to return a non-empty subset is called a \emph{choice function}.
The set of strategies that a choice function returns is called a \emph{normal form solution} of the decision tree.
So in the case of Fig.~\ref{fig:lake:tree:basic}, we would first find the six strategies and then apply our choice function to these.
Let us suppose that a particular choice function gives the normal form solution in Fig.~\ref{fig:normal:form:solution:example}.
Note that here, we apply the choice function conditional on the certain event.
Conditioning on other events is only relevant for subtrees, discussed next.

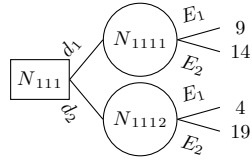
\begin{figure}
  \begin{center}
    \begin{tikzpicture}
      [minimum size=2em,parent anchor=east,child anchor=west,grow'=east,scale=0.75,transform shape]
            \node[draw,rectangle]{$\decnode[1](1)[1]$}
            [sibling distance=4em,level distance=5em]
            child{
              node[draw,circle]{$\chancenode[1](1)[1](1)$}
              [sibling distance=1.2em,level distance=5em]
              child{
                node{$9$}
                edge from parent
                node[above,sloped]{$E_1$}
              }
              child{
                node{$14$}
                edge from parent
                node[below,sloped]{$E_2$}
              }
              edge from parent
              node[above,sloped]{$d_1$}
            }
            child{
              node[draw,circle]{$\chancenode[1](1)[1](2)$}
              [sibling distance=1.2em,level distance=5em]
              child{
                node{$4$}
                edge from parent
                node[above,sloped]{$E_1$}
              }
              child{
                node{$19$}
                edge from parent
                node[below,sloped]{$E_2$}
              }
              edge from parent
              node[below,sloped]{$d_2$}
            }
            ;
    \end{tikzpicture}
  \end{center}
    \caption{The subtree of Fig.~\ref{fig:lake:tree:basic} at $\protect\decnode[1](1)[1]$.}
  \label{fig:subtree:example}
\end{figure}

\begin{figure}
  \begin{center}
    $\Bigg\{$
    \begin{tikzpicture}
      [minimum size=2em,parent anchor=east,child anchor=west,grow'=east,scale=0.75,transform shape,baseline=0em]
            \node[draw,rectangle]{$\decnode[1](1)[1]$}
            [sibling distance=4em,level distance=5em]
            child{
              node[draw,circle]{$\chancenode[1](1)[1](1)$}
              [sibling distance=1.2em,level distance=5em]
              child{
                node{$9$}
                edge from parent
                node[above,sloped]{$E_1$}
              }
              child{
                node{$14$}
                edge from parent
                node[below,sloped]{$E_2$}
              }
              edge from parent
              node[above,sloped]{$d_1$}
            }
            child{
              node{}
              edge from parent
              node[below,sloped]{$d_2$}
              node[sloped]{$\|$}
            }
            ;
    \end{tikzpicture}
    $\Bigg\}$
  \end{center}
    \caption{Possible normal form solution of Fig.~\ref{fig:subtree:example}.}
  \label{fig:strategy:example:subtree}
\end{figure}
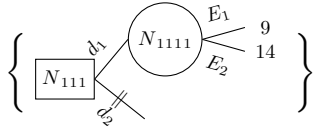

To examine subtree perfectness, we need to know how to restrict trees to smaller subtrees and solve these subtrees, and also how to restrict the \emph{solutions} of larger trees to smaller subtrees.
The first task is very straightforward: for each node on a tree, the subtree at that node is the node itself and all of its descendants.
For example, the subtree at $\decnode[1](1)[1]$ in Fig.~\ref{fig:lake:tree:basic} is shown in Fig.~\ref{fig:subtree:example}.
Observe that the labels of the nodes are as in the original tree.
Then the final step of the procedure is simply to solve this subtree,
conditional on the logical conjunction of all events leading up to $\decnode[1](1)[1]$ (or, conditional on the certain event, if no chance nodes are visited).\footnote{%
  We exclude decision trees that have $\emptyset$
  for the conditioning event of any of their subtrees,
  as most theories of conditioning do not allow conditioning on the empty set.
  This goes without loss of generality,
  because any decision tree can be `fixed'
  simply by pruning its $\emptyset$-conditioned subtrees.}
In our example, we find the two normal form solutions in Fig.~\ref{fig:subtree:example} and then apply our choice function, conditional on $S_1$ (this is the only event leading up to $\decnode[1](1)[1]$), to this set.
Let us suppose this gives the singleton strategy in Fig.~\ref{fig:strategy:example:subtree}.

This process of restriction and then solution can be summarized as follows.

\begin{enumerate}
  \item Choose a node $N$ to which to restrict the tree $\tree$.
  \item Find the conjunction of events which lead up to $N$
    from the root of $\tree$.
    Denote this conjunction by $E$.
    If there are no events leading up to $N$,
    then take $E$ to be the certain event.
  \item Remove from $\tree$ all non-descendants of $N$ other than $N$ itself. This is the subtree at $N$.
  \item Find all strategies for this subtree.
  \item Apply the choice function to this set of strategies, conditional on $E$.
    This is the normal form solution of the restricted tree.
\end{enumerate}

The opposite order, solving the large tree and then restricting the solution, is a little more involved.
We saw the first step in Fig.~\ref{fig:normal:form:solution:example}.
How can this solution be restricted to, say, $\decnode[1](1)[1]$?
This is done element-by-element.
We take the first strategy in the solution, and restrict it to $\decnode[1](1)[1]$ exactly as above.
This gives the same tree as in Fig.~\ref{fig:strategy:example:subtree}.
We then attempt to restrict the second strategy to $\decnode[1](1)[1]$, but we notice that this node does not appear in the strategy.
In such a case, we remove the strategy entirely.
Next we move to the third strategy, but again $\decnode[1](1)[1]$ is not present, so remove this strategy too.
This is the final strategy in the solution, so we find that the restriction of the normal form solution to $\decnode[1](1)[1]$ is the singleton in Fig.~\ref{fig:strategy:example:subtree}.
In this example, we see that restricting the tree to $\decnode[1](1)[1]$ and then solving is the same as solving the full tree and restricting the solution to $\decnode[1](1)[1]$.
We say that subtree perfectness is satisfied at $\decnode[1](1)[1]$.

In general, the process for solving a tree and then restricting the solution to a subtree goes as follows.
\begin{enumerate}
  \item Find all strategies for the tree.
  \item Apply the choice function to this set of strategies. This is the normal form solution of the tree.
  \item Choose a node $N$ to restrict the solution to.
  \item For each strategy $\atree$ in the solution:
  \begin{itemize}
    \item If $N$ is not a node of $\atree$, remove $\atree$.
    \item Otherwise, remove from $\atree$ all non-descendants of $N$ other than $N$ itself.
  \end{itemize}
  \item This gives the restriction of the normal form solution.
\end{enumerate}

Now we have all the tools to define subtree perfectness.

\begin{definition}[Subtree Perfectness]
  Let a decision tree $\tree$, a node $N$ in $\tree$, and a choice function be given.
  Carry out the operations of restricting the solution of $\tree$ to $N$, and solving the restriction of $\tree$ to $N$.
  The choice function is called subtree perfect for $\tree$ at $N$ if either
  \begin{itemize}
    \item the restriction of the solution is empty (that is, $N$ is not part of any strategy in the solution of $\tree$); or
    \item the restriction of the solution equals the solution of the restriction.
  \end{itemize}
  The choice function is called subtree perfect for $\tree$ if it is subtree perfect at all nodes $N$ in $\tree$.
  The choice function is called subtree perfect if it is subtree perfect for every $\tree$.
\end{definition}

Subtree perfect solutions are useful for several reasons.
When solving a decision problem, one can usually see it as being embedded in some larger decision problem that started in the past.
With a subtree perfect solution, the rest of this larger problem can be ignored.
A choice function that fails subtree perfectness may require consideration of events that did not happen or options that were refused in the past.
This is practically and perhaps intuitively unappealing.
Similarly, if the decision tree changes midway through the problem (for instance, new options that were not modeled become available), lack of subtree perfectness may require considering the whole problem again, whereas with subtree perfectness only the local problem needs to be solved.

\subsection{Choice Functions on Gambles}\label{sec:definition:choice:function}

The procedure introduced for finding a normal form solution involved applying a choice function to the set of available strategies.
In practice, it is rare to apply choice functions to sets of strategies, but rather to sets of \emph{gambles}.
A gamble is a map from the set of possible outcomes to a set of
rewards---to avoid technicalities,
we assume the set of possible outcomes to be finite.

For example,
consider the three strategies in Fig.~\ref{fig:normal:form:solution:example}.
The first of these will give $9$ if $S_1\cap E_1$ occurs, $14$ if $S_1\cap E_2$ occurs, $4$ if $S_2\cap E_1$ occurs, and $19$ if $S_2\cap E_2$ occurs.
This covers all possible outcomes, so for any outcome we know what reward the subject will receive.
This is exactly what we call a gamble.

We write such a gamble $9S_1E_1+14S_1E_2+4S_2E_1+19S_2E_2$, or even as $S_1(9E_1+14E_2)+S_2(4E_1+19E_2)$.
Similarly, the other two strategies in Fig.~\ref{fig:normal:form:solution:example} correspond to the gambles $S_1(10E_1+15E_2)+S_2(10E_1+15E_2)$ and $S_1(5E_1+20E_2)+S_2(5E_1+20E_2)$ respectively.
Note that which one of $S_1$ and $S_2$ occurs is irrelevant to both of these gambles, so they could be written as $10E_1+15E_2$ and $5E_1+20E_2$ respectively.
These objects are much more natural to work with; for instance, should we know the probabilities of the four events, then we could calculate the expected utility of the three gambles.
To work with sequential decision making, we actually need \emph{conditional} choice functions: the set returned by the choice function depends on what event is being conditioned upon.

Many, if not all, common choice functions used in decision theory are choice functions on gambles, hence restricting attention to them makes sense.
This adds two extra steps to the procedure for finding a normal form solution:
\begin{enumerate}
  \item Find all available strategies.
  \item For each strategy, find its corresponding gamble.
  \item Apply the choice function to this set of gambles.
  \item For each gamble chosen by the choice function, find all strategies in the original set that induced this gamble. This is the normal form solution. In other words, a strategy is in the normal form solution if and only if its corresponding gamble is chosen by the choice function.
\end{enumerate}
Note that we lose some generality by moving to choice functions on gambles, because there may be several strategies with the same corresponding gamble.
A choice function on gambles would be unable to distinguish between these, but a choice function on strategies would.
This loss is not so important, given that it would be rare to want to distinguish between strategies that give the same reward in every eventuality.

\subsection{Subtree Perfectness Theorem}\label{sec:theorem}

\begin{figure}
    \begin{center}
      \subfloat[]{
      \begin{tikzpicture}
        [minimum size=2em,parent anchor=east,child anchor=west,grow'=east,transform shape]
        \node[draw,circle]{$\chancenode(1)$}
        [level distance = 4em,sibling distance=3em]
          child{
            node[draw,rectangle]{$\decnode(1)[1]$}
            [sibling distance = 1.2em]
            child{
              node[right]{$X_1$}
            }
            child{
              node[right]{\tiny$\vdots$}
            }
            child{
              node[right]{$X_n$}
            }
            edge from parent
            node[above,sloped]{$A$}
          }
          child{
            node[right]{$Z$}
            edge from parent
            node[below,sloped]{$\compl{A}$}
          };        
      \end{tikzpicture}
      \label{fig:subtree:perfecness:simple:trees:1}
      }
      \hspace{1em}
      \subfloat[]{
      \begin{tikzpicture}
      [minimum size=2em,parent anchor=east,child anchor=west,grow'=east,transform shape]
      \node[draw,rectangle]{$\decnode[2]$}
      [level distance = 4em, sibling distance=4em]
        child{
          node[draw,rectangle]{$\decnode[2][1]$}
          [sibling distance = 1.2em]
          child{
            node[right]{$X_1$}
          }
          child{
            node[right]{\tiny$\vdots$}
          }
          child{
            node[right]{$X_m$}
          }
        }
        child{
          node[draw,rectangle]{$\decnode[2][2]$}
          [sibling distance = 1.2em]
          child{
            node[right]{$Y_1$}
          }
          child{
            node[right]{\tiny$\vdots$}
          }
          child{
            node[right]{$Y_n$}
          }
        };        
    \end{tikzpicture}
    \label{fig:subtree:perfecness:simple:trees:2}
    }
    \caption{Decision trees for which subtree perfectness implies subtree perfectness for any other tree.}
    \label{fig:subtree:perfecness:simple:trees}
  \end{center}
\end{figure}
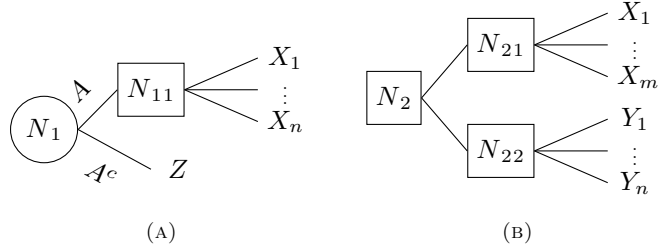

Now we can present our main result.
We are interested in when a particular choice function on gambles is subtree perfect.
It turns out that all that is required for this is subtree perfectness for two simple classes of decision trees, shown in Fig.~\ref{fig:subtree:perfecness:simple:trees}.
In that figure, the $X_i$, $Y_i$ and $Z$ are arbitrary gambles.
As part of a tree, they represent a terminal chance node.
For instance, the gamble $X=5E_1+3E_2$ would represent a chance node with two arcs,
with event $E_1$ leading to reward $5$, and event $E_2$ leading to reward $3$.

\begin{theorem}[Subtree Perfectness Theorem]
  \label{thm:subtreeperfectness}
  A choice function on gambles is subtree perfect if and only if it is subtree perfect for all decision trees of the two types in Fig.~\ref{fig:subtree:perfecness:simple:trees}, for every conditioning event $B$ at $\chancenode(1)$ such that $A\cap B\neq\emptyset$ and $\compl{A}\cap B\neq\emptyset$,
  and
  for every non-empty conditioning event $B$ at $\decnode[2]$.
\end{theorem}

So, in Fig.~\ref{fig:subtree:perfecness:simple:trees},
we consider solving the full trees conditional not just on the certain event,
but conditional on every possible non-empty event
that does not cause the decision tree to have a path
corresponding to $\emptyset$.

This theorem is equivalent to the subtree perfectness theorem presented and proved in \cite[Theorem~22]{2011:huntley:subtree:perfectness},
but is perhaps more straightforward to understand.
All we need to check for subtree perfectness is two simple types of decision tree, and this proves very useful when seeking both proofs and counterexamples for a particular choice function.
Note that it is not immediately obvious that this theorem is equivalent to \cite[Theorem~22]{2011:huntley:subtree:perfectness}, but in fact it appears in the proof \cite[A.9]{2011:huntley:subtree:perfectness}.

\section{Examples}\label{sec:examples}

\subsection{Examples of Subtree Perfectness}\label{sec:subtreeperfectness:examples}

In this section we give some examples of choice functions and decision trees where subtree perfectness holds.
We concentrate on two simple trees, shown in Fig.~\ref{fig:subtree:perfecness:simple:trees}.
This is not just for simplicity;
as we saw in Section~\ref{sec:subtreeperfectness},
these two types of trees are key to subtree perfectness.

\subsubsection{Expected Utility}\label{sec:subtreeperfectness:examples:eu}

\begin{example}\label{ex:eu:success:decnodes}
  First, we consider the simplest case, and investigate the tree
  in Fig.~\ref{fig:subtree:perfecness:simple:trees}\subref{fig:subtree:perfecness:simple:trees:2}, with conditioning event $B$.
  Maximizing conditional expectation for this tree is always subtree perfect,
  simply because
  \begin{multline}
    \max\{E(X_1|B),\dots,E(X_n|B),E(Y_1|B),\dots,E(Y_m|B)\}
    \\
    =
    \max\{\max\{E(X_1|B),\dots,E(X_n|B)\},\max\{E(Y_1|B),\dots,E(Y_m|B)\}\}.
  \end{multline}
\end{example}

\begin{example}\label{ex:eu:success:chancenodes}
  Consider the tree in Fig.~\ref{fig:subtree:perfecness:simple:trees}\subref{fig:subtree:perfecness:simple:trees:1}, with conditioning event $B$
  such that $A\cap B\neq\emptyset$ and $\compl{A}\cap B\neq\emptyset$.
  On the one hand, the normal form solution of the full tree is
  obtained by those gambles $AX_i+\compl{A}Z$ which maximize their
  expectation conditional on $B$:
  \begin{equation}
    \arg\max_{i=1}^nE(AX_i+\compl{A}Z|B).
  \end{equation}
  On the other hand,
  the normal form solution of the subtree at $\decnode(1)[1]$ is obtained by
  those gambles $X_i$ which maximize their expectation
  conditional on $A\cap B$:
  \begin{equation}
    \arg\max_{i=1}^nE(X_i|A\cap B)
  \end{equation}
  Now, note that, by the partition theorem,
  \begin{equation}
    E(AX_i+\compl{A}Z|B)=P(A|B)E(X_i|A\cap B)+P(\compl{A}|B)E(Z|\compl{A}\cap B).
  \end{equation}
  Because the term $P(\compl{A}|B)E(Z|\compl{A}\cap B)$ is independent of the index $i$,
  the solution of both maximization problems is identical:
  \begin{equation}
    \arg\max_{i=1}^nE(AX_i+\compl{A}Z|B)=\arg\max_{i=1}^nE(X_i|A\cap B),
  \end{equation}
  provided that $P(A|B)>0$. In other words, maximizing expected utility
  satisfies subtree perfectness on this type of tree
  whenever the probability of every non-empty event is strictly positive.
\end{example}

From the above analysis, and Theorem~\ref{thm:subtreeperfectness},
we recover a well-known result:
expected utility satisfies subtree perfectness for arbitrary decision trees,
provided that all chance arcs have strictly positive probability.
Interestingly, if some chance arcs have probability zero,
subtree perfectness may fail;
we already saw an instance of such failure
in the last example of Section~\ref{sec:introductory:example}.

\subsubsection{Choice Functions Induced By A Total Preorder}
\label{sec:subtreeperfectness:examples:total:preorder}

In Example~\ref{ex:eu:success:decnodes}, where we established
subtree perfectness of maximizing expectation for the tree
in Fig.~\ref{fig:subtree:perfecness:simple:trees}\subref{fig:subtree:perfecness:simple:trees:2},
all we really relied on is the fact that
expectation induces a total preorder over gambles.

Let us be a bit more specific.
Consider a relation $\succeq$ over gambles,
which satisfies transitivity, reflexivity, and completeness; that is,
for all gambles $X$, $Y$, and $Z$, we require that:
\begin{enumerate}
\item $X\succeq Y$ and $Y\succeq Z$ implies $X\succeq Z$;
\item $X\succeq X$; and
\item $X\succeq Y$ or $Y\succeq X$.
\end{enumerate}
Any relation which satisfies these properties is called a \emph{total preorder}.
Now, an obvious way to associate a choice function with a total preorder
is to select those gambles which dominate all other gambles:
\begin{equation}
  \max_{\succeq}\mathcal{X}\coloneqq \{X\in\mathcal{X}\colon\forall Y\in\mathcal{X}, X\succeq Y\},
\end{equation}
for any finite set of gambles $\mathcal{X}$.
Now, it is easily established that,
for any total preorder $\succeq$ over gambles,
\begin{equation}
  \max_{\succeq}\{X_1,\dots,X_n,Y_1,\dots,Y_m\}
  =
  \max_{\succeq}\{\max_{\succeq}\{X_1,\dots,X_n\},\max_{\succeq}\{Y_1,\dots,Y_m\}\}.
\end{equation}
In other words, any choice function induced by a total preorder
is subtree perfect for trees of the type of the tree
in Fig.~\ref{fig:subtree:perfecness:simple:trees}\subref{fig:subtree:perfecness:simple:trees:2}.

In fact, a much stronger result holds:
a choice function is subtree perfect \emph{only if}
it is induced by a total preorder
\cite[Theorem~22 \&\ Lemma~15]{2011:huntley:subtree:perfectness}.
Consequently, many of the choice functions used for
imprecise probability and imprecise utility fail subtree perfectness,
as they are typically not induced by total preorders; specifically,
they typically drop completeness.
Interestingly, completeness is not required
for some schemes of backward induction to work
\cite{2005:kikuti,2008:huntley:troffaes::impdectrees:smps,2011:kikuti:sequential,2012:huntley::backinduct}.

\subsubsection{Maximin}\label{sec:subtreeperfectness:examples:maximin}

The maximin criterion values a gamble at its smallest possible reward,
and then chooses gambles that give the maximum value.
Because it is a total preorder,
subtree perfectness will hold for the tree in Fig.~\ref{fig:subtree:perfecness:simple:trees}\subref{fig:subtree:perfecness:simple:trees:2},
as seen in Section~\ref{sec:subtreeperfectness:examples:total:preorder}.

Maximin does not in general hold for the tree in Fig.~\ref{fig:subtree:perfecness:simple:trees}\subref{fig:subtree:perfecness:simple:trees:1};
we shall see this later in Example~\ref{ex:maximin:failure:chancenodes}.

\subsubsection{$\Gamma$-Maximin}\label{sec:subtreeperfectness:examples:gamma:maximin}

Suppose that too little information is available
to identify a unique probability distribution, but we can specify a set of plausible probability distributions (see for instance \cite{1854:boole,1975:williams:condprev,1980:levi,2007:williams:condprev,1991:walley}).
In particular, suppose that
\newcommand{\pralt}{P}
\begin{itemize}
\item rewards are expressed in utiles, so $\rewardset=\SetR$,
\item the subject can express her beliefs by means of a set $\domlinprevs$ of probability distributions $\pralt$ ($\domlinprevs$ is called the \emph{credal set}, and is typically assumed to be closed and convex, however we need no such assumption here), and
\item each probability distribution $\pralt\in\domlinprevs$ satisfies $\pralt(A)>0$ for all events $A\neq\emptyset$.
\item  Note that the subject does not express beliefs about the relative plausibility of the distributions, so there is no second-order probability over the set. Also note that the positive probabilities are not necessary to define imprecise probability, but this assumption lets us avoid some technical details.
\end{itemize}

Under the above assumptions, each $\pralt$ in $\domlinprevs$ determines a conditional expectation\footnote{In this paper, gambles, such as $X$, are always assumed to take only a finite number of values.}
\begin{equation*}
  E_\pralt(X|A)=\frac{\sum_{x\in\SetR}x\pralt(X=x)}{\pralt(A)},
\end{equation*}
and the whole set $\domlinprevs$ determines a conditional lower and upper expectation
\begin{align}
  \lnex(X|A)&=\min_{\pralt\in\domlinprevs}E_\pralt(X|A)
  &
  \unex(X|A)&=\max_{\pralt\in\domlinprevs}E_\pralt(X|A), \label{eq:lower:envelope}
\end{align}
and this for every gamble $X$ and every non-empty event $A$.
We consider four choice functions
for this uncertainty model.
Each is introduced only briefly;
see \cite{2007:troffaes} for a more in-depth discussion and comparison.

The first, $\Gamma$-maximin, selects any gamble that maximizes lower expectation $\lnex$.
Failure of subtree perfectness for $\Gamma$-maximin has been well-documented~\cite{1988:seidenfeld, 2004:seidenfeld}, but for lower expectations with a particular structure, subtree perfectness can hold.

\begin{example}\label{ex:gammamaximin:success:chancenodes}
  Consider the tree in Fig.~\ref{fig:subtree:perfecness:simple:trees}\subref{fig:subtree:perfecness:simple:trees:1};
  for simplicity of notation, we assume here that the conditioning event $B$
  is the certain event, however the analysis extends trivially
  to any conditioning event.

  Suppose also that $\lnex$ satisfies
  \emph{marginal extension}~\cite[\S 6.7.2]{1991:walley}, that is,
  \begin{equation}\label{eq:marginal:extension}
    \lnex(A X + \compl{A} Y)=\lnex(A \lnex(X|A)+\compl{A}\lnex(Y|\compl{A})).
  \end{equation}
  for any gambles $X$ and $Y$.
  Note that this is naturally satisfied if beliefs are expressed as an
  unconditional lower expectation over gambles of the form $A\alpha+\compl{A}\beta$
  (for arbitrary $\alpha$, $\beta\in\SetR$), and
  conditional lower expectations $\lnex(\cdot|A)$ and
  $\lnex(\cdot|\compl{A})$.

  On the one hand, the normal form solution of the full tree is
  obtained by those gambles $AX_i+\compl{A}Z$ which maximize their lower expectation:
  \begin{equation}
    \arg\max_{i=1}^n\lnex(AX_i+\compl{A}Z).
  \end{equation}
  On the other hand,
  the normal form solution of the subtree at $\decnode(1)[1]$ is obtained by
  those gambles $X_i$ which maximize their lower expectation
  conditional on $A$:
  \begin{equation}
    \arg\max_{i=1}^n\lnex(X_i|A)
  \end{equation}
  Now, note that, by marginal extension,
  \begin{equation}
    \lnex(AX_i+\compl{A}Z)=\lnex(A \lnex(X_i|A)+\compl{A}\lnex(Z|\compl{A})).
  \end{equation}
  Using the above identity,
  one can easily show that the solution of both maximization problems is identical:
  \begin{equation}
    \arg\max_{i=1}^n\lnex(AX_i+\compl{A}Z)=\arg\max_{i=1}^n\lnex(X_i|A)
  \end{equation}
  provided that $\lnex(A)>0$. In other words, $\Gamma$-maximin
  satisfies subtree perfectness on this type of tree
  whenever the lower probability of every event is strictly positive.
\end{example}

Clearly, $\Gamma$-maximin corresponds to a total preorder, and therefore,
it will also satisfy subtree perfectness for the tree
in Fig.~\ref{fig:subtree:perfecness:simple:trees}\subref{fig:subtree:perfecness:simple:trees:2}.
Concluding,
as long as the lower probability of every event is strictly positive,
and marginal extension is satisfied with respect to the way the chance nodes partition the possibility space,
$\Gamma$-maximin satisfies subtree perfectness.

\subsubsection{$\Gamma$-Maximax}

The next criterion, $\Gamma$-maximax,
selects any gamble that maximizes upper expectation $\unex$.
The treatment of this criterion is very similar to that of $\Gamma$-maximax.
Here, we simply mention that $\Gamma$-maximax satisfies subtree perfectness
under exactly the same conditions as $\Gamma$-maximin.

\subsubsection{Maximality}

Maximality
starts from a strict partial order between gambles, conditional on an event $A$:
\begin{equation}
  X\succ_A Y\text{ whenever }\lnex(X-Y|A)>0
\end{equation}
and then selects those gambles which are undominated with respect to this strict partial order:
\begin{equation}
  \max_{\succ_A}\mathcal{X}
  \coloneqq
  \{X\in\mathcal{X}\colon \forall Y\in\mathcal{X},\,Y\not\succ_A X\}
\end{equation}
for any finite set of gambles $\mathcal{X}$.

Here, one can show that,
under strictly positive lower probability of conditioning events,
maximality will satisfy subtree perfectness
for trees of the type of Fig.~\ref{fig:subtree:perfecness:simple:trees}\subref{fig:subtree:perfecness:simple:trees:1}, essentially because
\begin{equation}
  \lnex\big((AX_i+\compl{A}Z)-(AX_j+\compl{A}Z)\big| B\big)
  =\lnex\big(A (X_i-X_j)\big| B\big).
\end{equation}
has the same sign as $\lnex(X_i-X_j|A\cap B)$, because of the generalized Bayes rule
\cite{1991:walley},
provided that $\lnex(A|B)>0$. Note that we do not require marginal extension.

We will see later that it does not satisfy subtree perfectness for trees of the type of Fig.~\ref{fig:subtree:perfecness:simple:trees}\subref{fig:subtree:perfecness:simple:trees:2}.

In other words, maximality satisfies subtree perfectness
only in some special cases.
Interestingly,
there are backward induction schemes which do work for maximality
(under strictly positive lower probability)
\cite{2005:kikuti,2008:huntley:troffaes::impdectrees:smps,2011:kikuti:sequential,2012:huntley::backinduct}.

\subsubsection{E-admissibility}\label{sec:subtreeperfectness:examples:E:admissibility}

E-admissibility, selects any gamble
that maximizes expectation
under at least one probability distribution in $\domlinprevs$.

Under strictly positive lower probability,
as with maximality,
it satisfies  subtree perfectness for trees of the type of Fig.~\ref{fig:subtree:perfecness:simple:trees}\subref{fig:subtree:perfecness:simple:trees:1}.
It does not satisfy subtree perfectness for trees of the type of Fig.~\ref{fig:subtree:perfecness:simple:trees}\subref{fig:subtree:perfecness:simple:trees:2}.

The next example uses E-admissibility to demonstrate that subtree perfectness
is not as strong a condition as it may first appear.

\begin{figure}
  \begin{center}
    \begin{tikzpicture}
      [minimum size=2em,parent anchor=east,child anchor=west,grow'=east,transform shape]
      \node[draw,circle]{$N$}
      [level distance = 4em,sibling distance=8em]
      child{
        node[draw,rectangle]{$N_1$}
        [sibling distance = 3em]
        child{
          node[draw,circle]{$N_{11}$}
          child{
            node[right]{$1$}
            edge from parent
            node[sloped,above]{$B_1$}
          }
          child{
            node[right]{$-1$}
            edge from parent
            node[sloped,below]{$B_2$}
          }
        }
        child{
          node[right]{$0$}
        }
        edge from parent
        node[above,sloped]{$A_1$}
      }
      child{
        node[draw,rectangle]{$N_2$}
        [sibling distance = 3em]
        child{
          node[draw,circle]{$N_{21}$}
          child{
            node[right]{$-1$}
            edge from parent
            node[sloped,above]{$B_1$}
          }
          child{
            node[right]{$1$}
            edge from parent
            node[sloped,below]{$B_2$}
          }
        }
        child{
          node[right]{$0$}
        }
        edge from parent
        node[below,sloped]{$A_2$}
      };
    \end{tikzpicture}
    \caption{A tree to demonstrate that subtree perfectness is not as strong as it may first appear.}
    \label{fig:subtree:perfecness:eadm:tree}
  \end{center}
\end{figure}
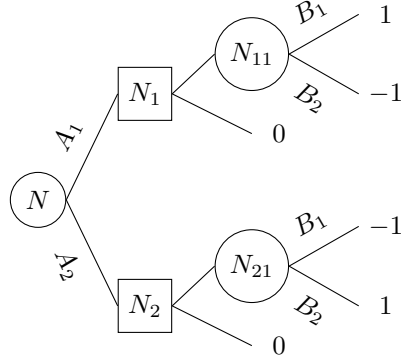

\begin{figure}
  \begin{center}
    $\left\{\rule{0em}{6em}\right.$
    \begin{tikzpicture}
      [minimum size=2em,parent anchor=east,child anchor=west,grow'=east,transform shape,baseline=0em]
      \node[draw,circle]{$N$}
      [level distance = 4em,sibling distance=8em]
      child{
        node[draw,rectangle]{$N_1$}
        [sibling distance = 3em]
        child{
          node[draw,circle]{$N_{11}$}
          child{
            node[right]{$1$}
            edge from parent
            node[sloped,above]{$B_1$}
          }
          child{
            node[right]{$-1$}
            edge from parent
            node[sloped,below]{$B_2$}
          }
        }
        child{
          node[right]{}
          edge from parent
          node[sloped]{$\|$}
        }
        edge from parent
        node[above,sloped]{$A_1$}
      }
      child{
        node[draw,rectangle]{$N_2$}
        [sibling distance = 3em]
        child{
          node[right]{}
          edge from parent
          node[sloped]{$\|$}
        }
        child{
          node[right]{$0$}
        }
        edge from parent
        node[below,sloped]{$A_2$}
      };
    \end{tikzpicture}
    \hfill,\hfill
    \begin{tikzpicture}
      [minimum size=2em,parent anchor=east,child anchor=west,grow'=east,transform shape,baseline=0em]
      \node[draw,circle]{$N$}
      [level distance = 4em,sibling distance=8em]
      child{
        node[draw,rectangle]{$N_1$}
        [sibling distance = 3em]
        child{
          node{}
          edge from parent
          node[sloped]{$\|$}
        }
        child{
          node[right]{$0$}
        }
        edge from parent
        node[above,sloped]{$A_1$}
      }
      child{
        node[draw,rectangle]{$N_2$}
        [sibling distance = 3em]
        child{
          node[draw,circle]{$N_{21}$}
          child{
            node[right]{$-1$}
            edge from parent
            node[sloped,above]{$B_1$}
          }
          child{
            node[right]{$1$}
            edge from parent
            node[sloped,below]{$B_2$}
          }
        }
        child{
          node{}
          edge from parent
          node[sloped]{$\|$}
        }
        edge from parent
        node[below,sloped]{$A_2$}
      };
    \end{tikzpicture}
    $\left.\rule{0em}{6em}\right\}$
    \caption{Solution of the tree of Fig.~\ref{fig:subtree:perfecness:eadm:tree}.}
    \label{fig:subtree:perfecness:eadm:tree:solution}
  \end{center}
\end{figure}
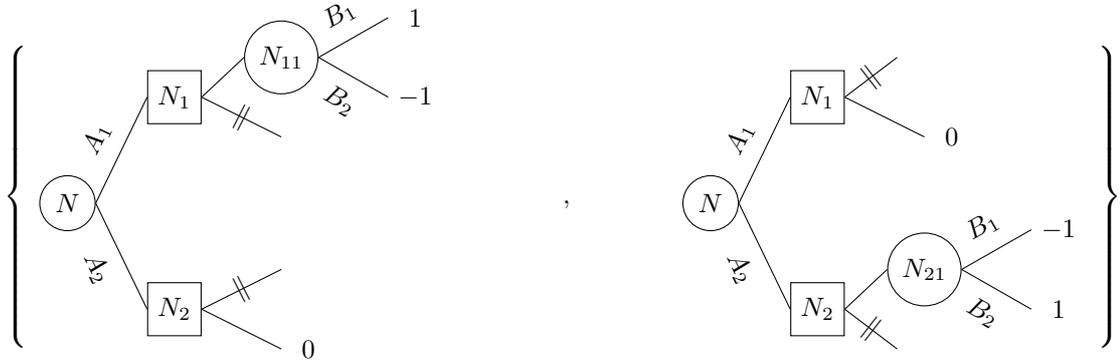

\begin{figure}
  \begin{center}
    \subfloat[]{
      $\Bigg\{$
      \begin{tikzpicture}
        [minimum size=2em,parent anchor=east,child anchor=west,grow'=east,transform shape,baseline=0em]
        \node[draw,rectangle]{$N_1$}
        [sibling distance = 3em]
        child{
          node[draw,circle]{$N_{11}$}
          child{
            node[right]{$1$}
            edge from parent
            node[sloped,above]{$B_1$}
          }
          child{
            node[right]{$-1$}
            edge from parent
            node[sloped,below]{$B_2$}
          }
        }
        child{
          node{}
          edge from parent
          node[sloped]{$\|$}
        }
        ;
      \end{tikzpicture}
      ,
      \begin{tikzpicture}
        [minimum size=2em,parent anchor=east,child anchor=west,grow'=east,transform shape,baseline=0em]
        \node[draw,rectangle]{$N_1$}
        [sibling distance = 3em]
        child{
          node{}
          edge from parent
          node[sloped]{$\|$}
        }
        child{
          node{$0$}
        }
        ;
      \end{tikzpicture}
      $\Bigg\}$
    }
    \\
    \subfloat[]{
      $\Bigg\{$
      \begin{tikzpicture}
        [minimum size=2em,parent anchor=east,child anchor=west,grow'=east,transform shape,baseline=0em]
        \node[draw,rectangle]{$N_2$}
        [sibling distance = 3em]
        child{
          node[draw,circle]{$N_{21}$}
          child{
            node[right]{$-1$}
            edge from parent
            node[sloped,above]{$B_1$}
          }
          child{
            node[right]{$1$}
            edge from parent
            node[sloped,below]{$B_2$}
          }
        }
        child{
          node{}
          edge from parent
          node[sloped]{$\|$}
        }
        ;
      \end{tikzpicture}
      ,
      \begin{tikzpicture}
        [minimum size=2em,parent anchor=east,child anchor=west,grow'=east,transform shape,baseline=0em]
        \node[draw,rectangle]{$N_2$}
        [sibling distance = 3em]
        child{
          node{}
          edge from parent
          node[sloped]{$\|$}
        }
        child{
          node{$0$}
        }
        ;
      \end{tikzpicture}
      $\Bigg\}$
    }
    \caption{Solutions of the subtrees at $N_1$ and $N_2$, for the tree of Fig.~\ref{fig:subtree:perfecness:eadm:tree}.}
    \label{fig:subtree:perfecness:eadm:tree:solution:subtrees}
  \end{center}
\end{figure}
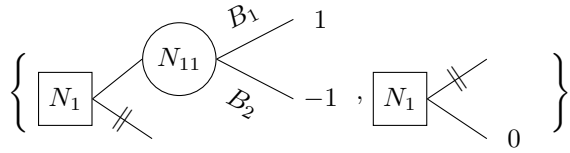
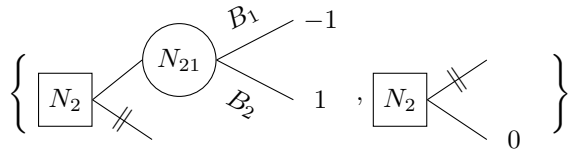

\begin{example}\label{ex:eadmissibility:success:special}
  Consider the tree in Fig.~\ref{fig:subtree:perfecness:eadm:tree}.
  Suppose that $\domlinprevs=\{\pr_1,\pr_2\}$, with
  $\pr_1(A_1)=\pr_1(A_2)=\pr_2(A_1)=\pr_2(A_2)=0.5$,
  $\pr_1(B_1|A_1)=\pr_1(B_1|A_2)=0.6$, and
  $\pr_2(B_1|A_1)=\pr_2(B_1|A_2)=0.4$.

  The normal form solution of the full tree is depicted in
  Fig.~\ref{fig:subtree:perfecness:eadm:tree:solution}.
  The left-hand solution maximizes expectation under $\pr_1$,
  whilst the right-hand solution maximizes expectation under $\pr_2$.
  Of the four strategies, two are optimal according to E-admissibility.

  The normal form solutions of the subtrees at $N_1$ and $N_2$ are
  depicted in Fig.~\ref{fig:subtree:perfecness:eadm:tree:solution:subtrees}.
  In both subtrees, there are two possible strategies,
  and both are optimal according to E-admissibility.

  Subtree perfectness holds here, but observe that only two of the four global strategies are optimal, whilst all of the local strategies are optimal.
\end{example}

This shows that subtree perfectness is weaker than
what one might call
\emph{normal-extensive form equivalence}.
For further discussion of the relationship between
subtree perfectness and normal-extensive form equivalence,
we refer to \cite{2011:huntley:subtree:perfectness}.

\subsection{Failures of Subtree Perfectness}\label{sec:failures:subtreeperfectness:examples}

In Section~\ref{sec:introductory:example},
we already visited one well-known instance subtree perfectness failing, namely
expected utility, when events of zero probability are involved.
In this section, we provide many more examples.

\subsubsection{Imprecise Utility}
\label{sec:imprecise:utility}

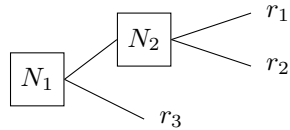
\begin{figure}
  \begin{center}
    \begin{tikzpicture}
      [minimum size=2em,parent anchor=east,child anchor=west,grow'=east,transform shape]
      \node[draw,rectangle]{$N_1$}
      [sibling distance=3em, level distance=4em]
      child{
        node[draw,rectangle]{$N_2$}
        [sibling distance=2em]
        child{
          node[right]{$r_1$}
        }
        child{
          node[right]{$r_2$}
        }
      }
      child{
        node[right]{$r_3$}
      };
    \end{tikzpicture}
    \caption{An example involving imprecise utility.}
    \label{fig:example:imprecise:utility}
  \end{center}
\end{figure}

\begin{example}\label{ex:imprecise:utility}
Consider again the very first example
from the introduction (Fig.~\ref{fig:two:stage:problem}),
redrawn in Fig.~\ref{fig:example:imprecise:utility};
note that this tree is an extremely simple instance of the tree in
Fig.~\ref{fig:subtree:perfecness:simple:trees}\subref{fig:subtree:perfecness:simple:trees:2}
which plays a central role in our subtree perfectness theorem.

Suppose we are unsure about our utility between the three available options,
and we are happy to accept any option as optimal that maximizes
either of the following two utility functions:

\begin{align}
  U_1(r_1)&=3,
  &
  U_1(r_2)&=1,
  &
  U_1(r_3)&=2,
  \\
  U_2(r_1)&=-3,
  &
  U_2(r_2)&=1,
  &
  U_2(r_3)&=2.
\end{align}
Clearly, for the full tree, the optimal choice is
$\{r_1,r_3\}$. However, for the subtree at $N_2$, the optimal choice
is $\{r_1,r_2\}$.  Subtree perfectness is violated, because $r_2$ is
not optimal in the full tree, whereas it is optimal for the subtree at $N_2$.
\end{example}

\subsubsection{Choice Functions Not Induced By A Total Preorder}\label{sec:subtree:perfectness:failure:not:total:preorder}

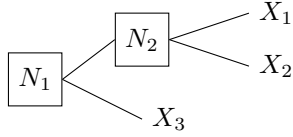
\begin{figure}
  \begin{center}
    \begin{tikzpicture}
      [minimum size=2em,parent anchor=east,child anchor=west,grow'=east,transform shape]
      \node[draw,rectangle]{$N_1$}
      [sibling distance=3em, level distance=4em]
      child{
        node[draw,rectangle]{$N_2$}
        [sibling distance=2em]
        child{
          node[right]{$X_1$}
        }
        child{
          node[right]{$X_2$}
        }
      }
      child{
        node[right]{$X_3$}
      };
    \end{tikzpicture}
    \caption{An example involving a partial order.}
    \label{fig:example:partial:order}
  \end{center}
\end{figure}

The previous example easily adapts to choice functions induced by a
strict partial order.
Indeed, assume that we have a strict partial order $\succ$ on gambles,
and three gambles, $X_1$, $X_2$, and $X_3$, such that:
\begin{equation}
  X_3\succ X_2,
\end{equation}
but where all other pairs are incomparable:
\begin{align}
  X_1&\not\succ X_2,
  &
  X_1&\not\succ X_3,
  \\
  X_2&\not\succ X_1,
  &
  X_3&\not\succ X_1.
\end{align}
Considering now the situation depicted in
Fig.~\ref{fig:example:partial:order},
we are in exactly the same situation as before:
clearly, for the full tree, the optimal choice is
$\{X_1,X_3\}$. However, for the subtree at $N_2$, the optimal choice
is $\{X_1,X_2\}$.  Subtree perfectness is violated, because $X_2$ is
dominated, and whence, not optimal, in the full tree,
whereas it is optimal for the subtree at $N_2$.

Although this is not a rigorous proof, it does
demonstrate why any non-trivial partial order might fail subtree perfectness.
In other words, maximality
(and also point-wise dominance and interval dominance,
which we have not discussed here, but which are also induced by partial orders
and are commonly used in the literature \cite{2007:troffaes})
will fail subtree perfectness in general.

Although E-admissibility does not correspond to a partial order,
unsurprisingly, it fails subtree perfectness in a very similar way.

\begin{figure}
  \begin{center}
    \begin{tikzpicture}
      [minimum size=2em,parent anchor=east,child anchor=west,grow'=east,transform shape]
      \node[draw,rectangle]{$N_1$}
      [sibling distance=3em, level distance=4em]
      child{
        node[draw,rectangle]{$N_2$}
        [sibling distance=3em]
        child{
          node[draw,circle]{$N_3$}
          [sibling distance=2em]
          child{
            node[right]{$5$}
            edge from parent
            node[sloped,above]{$A_1$}
          }
          child{
            node[right]{$-5$}
            edge from parent
            node[sloped,below]{$A_2$}
          }
        }
        child{
          node[right]{$1$}
        }
      }
      child{
        node[right]{$2$}
      };
    \end{tikzpicture}
    \caption{Failure of subtree perfectness for E-admissibility.}
    \label{fig:example:failure:eadm}
  \end{center}
\end{figure}
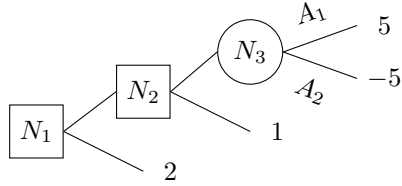

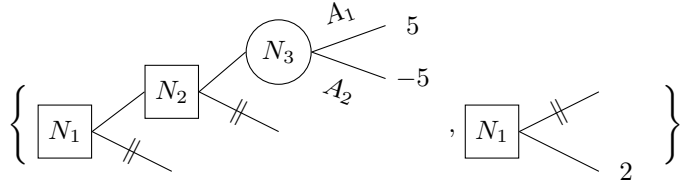
\begin{figure}
  \begin{center}
    $\Bigg\{$
    \begin{tikzpicture}
      [minimum size=2em,parent anchor=east,child anchor=west,grow'=east,transform shape,baseline=0em]
      \node[draw,rectangle]{$N_1$}
      [sibling distance=3em, level distance=4em]
      child{
        node[draw,rectangle]{$N_2$}
        [sibling distance=3em]
        child{
          node[draw,circle]{$N_3$}
          [sibling distance=2em]
          child{
            node[right]{$5$}
            edge from parent
            node[sloped,above]{$A_1$}
          }
          child{
            node[right]{$-5$}
            edge from parent
            node[sloped,below]{$A_2$}
          }
        }
        child{
          node[right]{}
          edge from parent
          node[sloped]{$\|$}
        }
      }
      child{
        node[right]{}
        edge from parent
        node[sloped]{$\|$}
      };
    \end{tikzpicture}
    ,
    \begin{tikzpicture}
      [minimum size=2em,parent anchor=east,child anchor=west,grow'=east,transform shape,baseline=0em]
      \node[draw,rectangle]{$N_1$}
      [sibling distance=3em, level distance=4em]
      child{
        node[right]{}
        edge from parent
        node[sloped]{$\|$}
      }
      child{
        node[right]{$2$}
      };
    \end{tikzpicture}
    $\Bigg\}$
    \caption{E-admissible normal form solution of the tree of Fig.~\ref{fig:example:failure:eadm}.}
    \label{fig:example:failure:eadm:solution}
  \end{center}
\end{figure}

\begin{figure}
  \begin{center}
    $\Bigg\{$
    \begin{tikzpicture}
      [minimum size=2em,parent anchor=east,child anchor=west,grow'=east,transform shape,baseline=0em]
        \node[draw,rectangle]{$N_2$}
        [sibling distance=3em]
        child{
          node[draw,circle]{$N_3$}
          [sibling distance=2em]
          child{
            node[right]{$5$}
            edge from parent
            node[sloped,above]{$A_1$}
          }
          child{
            node[right]{$-5$}
            edge from parent
            node[sloped,below]{$A_2$}
          }
        }
        child{
          node[right]{}
          edge from parent
          node[sloped]{$\|$}
        }
        ;
    \end{tikzpicture}
    ,
    \begin{tikzpicture}
      [minimum size=2em,parent anchor=east,child anchor=west,grow'=east,transform shape,baseline=0em]
      \node[draw,rectangle]{$N_2$}
      [sibling distance=3em, level distance=4em]
      child{
        node[right]{}
        edge from parent
        node[sloped]{$\|$}
      }
      child{
        node[right]{$1$}
      };
    \end{tikzpicture}
    $\Bigg\}$
    \caption{E-admissible normal form solution of the subtree at $N_2$, for the tree of Fig.~\ref{fig:example:failure:eadm}.}
    \label{fig:example:failure:eadm:solution:subtree}
  \end{center}
\end{figure}
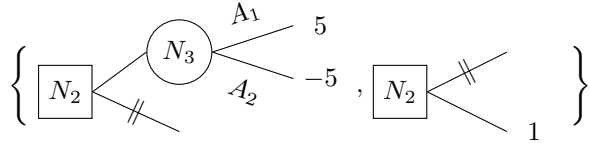

\begin{example}\label{ex:edmissibility:failure:decnodes}
	Consider the decision tree in Fig.~\ref{fig:example:failure:eadm}.
        Suppose that $\domlinprevs=\{\pr_1,\pr_2\}$,
        where $\pr_1(A_1)=\pr_2(A_2)=0.8$.
        The expectation of $N_3$ is $3$ under $\pr_1$
        and $-3$ under $\pr_2$. The remaining two strategies have
        constant values, namely $1$, and $2$.

        In the full tree, clearly, of the three strategies, only two
        are E-admissible,
        as depicted in Fig.~\ref{fig:example:failure:eadm:solution}.
        However, for the subtree at $N_2$, both strategies are
        optimal,
        as depicted Fig.~\ref{fig:example:failure:eadm:solution:subtree}.
        Subtree perfectness is violated, because
        the strategy with constant outcome $1$ is
        dominated by the strategy with constant outcome $2$,
        and whence, $1$ is not optimal in the full tree,
        whereas it is optimal for the subtree at $N_2$.
        In fact, we have mimicked the situation of
        Example~\ref{ex:imprecise:utility}.
\end{example}

\subsubsection{Maximin}
\label{sec:failures:subtreeperfectness:examples:maximin}

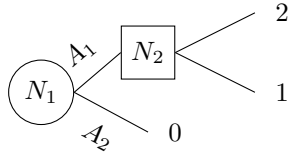
\begin{figure}
  \begin{center}
    \begin{tikzpicture}
      [minimum size=2em,parent anchor=east,child anchor=west,grow'=east,transform shape]
      \node[draw,circle]{$N_1$}
      [sibling distance=3em, level distance=4em]
      child{
        node[draw,rectangle]{$N_2$}
        [sibling distance=3em]
        child{
          node[right]{$2$}
        }
        child{
          node[right]{$1$}
        }
        edge from parent
        node[above,sloped]{$A_1$}
      }
      child{
        node[right]{$0$}
        edge from parent
        node[below,sloped]{$A_2$}
      };
    \end{tikzpicture}
    \caption{Failure of subtree perfectness for maximin.}
    \label{fig:example:maximin}
  \end{center}
\end{figure}

\begin{example}\label{ex:maximin:failure:chancenodes}
  Consider the tree in Fig.~\ref{fig:example:maximin};
  note that this tree is an extremely simple instance of the tree in
  Fig.~\ref{fig:subtree:perfecness:simple:trees}\subref{fig:subtree:perfecness:simple:trees:1}
  which plays a central role in our subtree perfectness theorem.

  Obviously, there are only two strategies. In the full tree, the worst possible
  outcome for both strategies is $0$, and therefore both strategies are maximin.
  However, clearly, at $N_2$, only the strategy yielding $2$ is optimal.
  Whence, subtree perfectness is violated.
\end{example}

\subsubsection{$\Gamma$-maximin}

When Eq.~\eqref{eq:marginal:extension} does not hold,
$\Gamma$-maximin often fails to be subtree perfect, in a similar way
that maximin fails it.

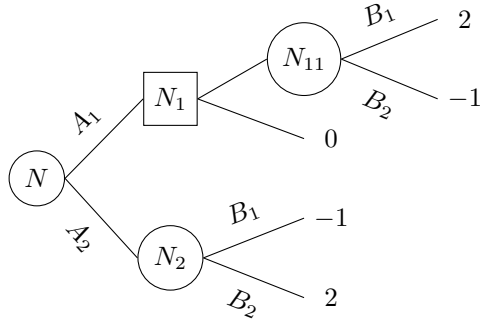
\begin{figure}
  \begin{center}
    \begin{tikzpicture}
      [minimum size=2em,parent anchor=east,child anchor=west,grow'=east,transform shape]
      \node[draw,circle]{$N$}
      [sibling distance=6em, level distance=5em]
      child{
        node[draw,rectangle]{$N_1$}
        [sibling distance=3em]
        child{
          node[draw,circle]{$N_{11}$}
          child{
            node[right]{$2$}
            edge from parent
            node[above,sloped]{$B_1$}
          }
          child{
            node[right]{$-1$}
            edge from parent
            node[below,sloped]{$B_2$}
          }
        }
        child{
          node[right]{$0$}
        }
        edge from parent
        node[above,sloped]{$A_1$}
      }
      child{
        node[draw,circle]{$N_2$}
        [sibling distance=3em]
        child{
          node[right]{$-1$}
          edge from parent
          node[above,sloped]{$B_1$}
        }
        child{
          node[right]{$2$}
          edge from parent
          node[below,sloped]{$B_2$}
        }
        edge from parent
        node[below,sloped]{$A_2$}
      };
    \end{tikzpicture}
    \caption{Failure of subtree perfectness for $\Gamma$-maximin.}
    \label{fig:example:gamma:maximin}
  \end{center}
\end{figure}

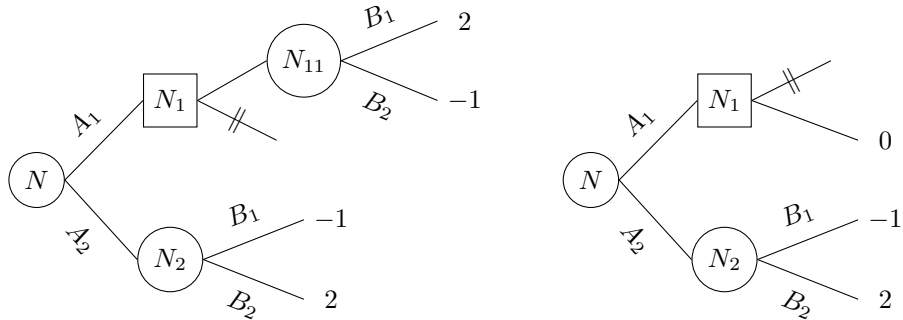
\begin{figure}
  \begin{center}
    \begin{tikzpicture}
      [minimum size=2em,parent anchor=east,child anchor=west,grow'=east,transform shape,baseline=0em]
      \node[draw,circle]{$N$}
      [sibling distance=6em, level distance=5em]
      child{
        node[draw,rectangle]{$N_1$}
        [sibling distance=3em]
        child{
          node[draw,circle]{$N_{11}$}
          child{
            node[right]{$2$}
            edge from parent
            node[above,sloped]{$B_1$}
          }
          child{
            node[right]{$-1$}
            edge from parent
            node[below,sloped]{$B_2$}
          }
        }
        child{
          node{}
          edge from parent
          node[sloped]{$\|$}
        }
        edge from parent
        node[above,sloped]{$A_1$}
      }
      child{
        node[draw,circle]{$N_2$}
        [sibling distance=3em]
        child{
          node[right]{$-1$}
          edge from parent
          node[above,sloped]{$B_1$}
        }
        child{
          node[right]{$2$}
          edge from parent
          node[below,sloped]{$B_2$}
        }
        edge from parent
        node[below,sloped]{$A_2$}
      };
    \end{tikzpicture}
    \hspace{2em}
    \begin{tikzpicture}
      [minimum size=2em,parent anchor=east,child anchor=west,grow'=east,transform shape,baseline=0em]
      \node[draw,circle]{$N$}
      [sibling distance=6em, level distance=5em]
      child{
        node[draw,rectangle]{$N_1$}
        [sibling distance=3em]
        child{
          node{}
          edge from parent
          node[sloped]{$\|$}
        }
        child{
          node[right]{$0$}
        }
        edge from parent
        node[above,sloped]{$A_1$}
      }
      child{
        node[draw,circle]{$N_2$}
        [sibling distance=3em]
        child{
          node[right]{$-1$}
          edge from parent
          node[above,sloped]{$B_1$}
        }
        child{
          node[right]{$2$}
          edge from parent
          node[below,sloped]{$B_2$}
        }
        edge from parent
        node[below,sloped]{$A_2$}
      };
    \end{tikzpicture}
    \caption{The two strategies for the tree of Fig.~\ref{fig:example:gamma:maximin}. Only the left-hand strategy is $\Gamma$-maximin optimal.}
    \label{fig:example:gamma:maximin:strategies}
  \end{center}
\end{figure}

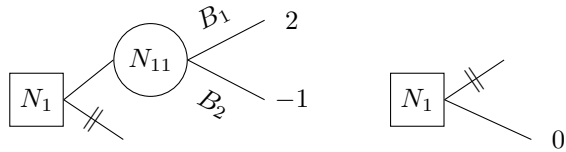
\begin{figure}
  \begin{center}
    \begin{tikzpicture}
      [minimum size=2em,parent anchor=east,child anchor=west,grow'=east,transform shape,baseline=0em]
        \node[draw,rectangle]{$N_1$}
        [sibling distance=3em]
        child{
          node[draw,circle]{$N_{11}$}
          child{
            node[right]{$2$}
            edge from parent
            node[above,sloped]{$B_1$}
          }
          child{
            node[right]{$-1$}
            edge from parent
            node[below,sloped]{$B_2$}
          }
        }
        child{
          node{}
          edge from parent
          node[sloped]{$\|$}
        }
        ;
    \end{tikzpicture}
    \hspace{2em}
    \begin{tikzpicture}
      [minimum size=2em,parent anchor=east,child anchor=west,grow'=east,transform shape,baseline=0em]
        \node[draw,rectangle]{$N_1$}
        [sibling distance=3em]
        child{
          node{}
          edge from parent
          node[sloped]{$\|$}
        }
        child{
          node[right]{$0$}
        }
        ;
    \end{tikzpicture}
    \caption{The two strategies for the subtree at $N_1$, for the tree of Fig.~\ref{fig:example:gamma:maximin}. Only the right-hand strategy is $\Gamma$-maximin optimal.}
    \label{fig:example:gamma:maximin:subtree:strategies}
  \end{center}
\end{figure}


\begin{example}\label{ex:gammamaximin:failure:chancenodes}
  Consider the tree in Fig.~\ref{fig:example:gamma:maximin};
  again
  note that this tree is a simple instance of the tree in
  Fig.~\ref{fig:subtree:perfecness:simple:trees}\subref{fig:subtree:perfecness:simple:trees:1}
  which plays a central role in our subtree perfectness theorem.

  Suppose that $\domlinprevs=\{\pr_1,\pr_2\}$ with
  \begin{center}
    \begin{tabular}{c|cc}
      $\pr_1$ & $A_1$ & $A_2$ \\
      \hline
      $B_1$ & $0.1$ & $0.1$ \\
      $B_2$ & $0.4$ & $0.4$
    \end{tabular}
    and
    \begin{tabular}{c|cc}
      $\pr_2$ & $A_1$ & $A_2$ \\
      \hline
      $B_1$ & $0.4$ & $0.4$ \\
      $B_2$ & $0.1$ & $0.1$
    \end{tabular}
  \end{center}

  The two strategies of the full tree are depicted in
  Fig.~\ref{fig:example:gamma:maximin:strategies},
  and correspond to the following gambles:
  \begin{center}
    \begin{tabular}{c|cc}
      & $A_1$ & $A_2$ \\
      \hline
      $B_1$ & $2$ & $-1$ \\
      $B_2$ & $-1$ & $2$
    \end{tabular}
    and
    \begin{tabular}{c|cc}
      & $A_1$ & $A_2$ \\
      \hline
      $B_1$ & $0$ & $-1$ \\
      $B_2$ & $0$ & $2$
    \end{tabular}
  \end{center}
  The lower expectation of the first strategy is $0.5$, 
  whereas the second strategy has lower expectation $-0.2$.
  So, in the full tree, only the first strategy is optimal
  according to $\Gamma$-maximin.

  However, in the subtree at $N_1$, the conditional probabilities are
  \begin{center}
    \begin{tabular}{c|cc}
      & $\pr_1(\cdot|A_1)$ \\
      \hline
      $B_1$ & $0.2$ \\
      $B_2$ & $0.8$
    \end{tabular}
    and
    \begin{tabular}{c|cc}
      & $\pr_2(\cdot|A_1)$ \\
      \hline
      $B_1$ & $0.8$ \\
      $B_2$ & $0.2$
    \end{tabular}
  \end{center}
  and two strategies, depicted in
  Fig.~\ref{fig:example:gamma:maximin:subtree:strategies},
  of this subtree, correspond to the following gambles:
  \begin{center}
    \begin{tabular}{c|cc}
      $B_1$ & $2$ \\
      $B_2$ & $-1$
    \end{tabular}
    and
    \begin{tabular}{c|cc}
      $B_1$ & $0$ \\
      $B_2$ & $0$
    \end{tabular}
  \end{center}
  The lower expectation of the first strategy is $-0.4$, 
  whereas the second strategy has lower expectation $0$.
  So, in the subtree, only the second strategy is optimal
  according to $\Gamma$-maximin.

  Concluding, $\Gamma$-maximin does not satisfy subtree perfectness.
  In fact, the preference between the two strategies reverses
  depending on whether we consider the full tree, or only the subtree at $N_1$.
\end{example}

We refer to Seidenfeld~\cite{2004:seidenfeld} for further beautiful examples.

\section{Conclusion}\label{sec:conclusion}

We introduced the concept of subtree perfectness: in essence, we extended Selten's idea of subgame perfectness to resolute normal form solutions of decision trees.
Subtree perfectness allows us to solve decision trees without worrying about which larger tree they are embedded in.
This is a particularly useful property 
if one desires to avoid counterfactual reasoning for philosophical reasons.
It also allows for efficient computational algorithms to solve decision trees
via backward induction,
and is tremendously helpful when something unexpected happens
in the middle of a sequential decision problem,
as subtree perfectness allows us to
restrict our investigation to a subtree of the larger problem.

We saw that, when using choice functions on gambles, subtree perfectness for all decision trees is equivalent to subtree perfectness for two simple classes of decision trees.
Numerous examples for popular choice functions, demonstrated how this result can be leveraged both for proving subtree perfectness and for finding simple counterexamples.
The application of the results to counterexamples is particularly useful, since it shows that, when searching for a counterexample, there is no reason to try complicated trees because a small one will always suffice.

For choice functions that are not total preorders (that is, those that attempt to model indeterminacy or indecision), we highlight a few important consequences.
First, such choice functions always fail subtree perfectness in general; see the examples in Sections~\ref{sec:imprecise:utility} and~\ref{sec:subtree:perfectness:failure:not:total:preorder} and~\cite[Lemmas~15 \& 21]{2011:huntley:subtree:perfectness} for a full proof.
We demonstrated that these counterexamples typically appear when the tree contains paths with more than one decision node.
However, and this may come as a surprise, such choice functions can satisfy subtree perfectness for decision trees where a chance node is followed by single decision nodes; see the examples in Sections~\ref{sec:subtreeperfectness:examples:gamma:maximin}--\ref{sec:subtreeperfectness:examples:E:admissibility}.
Most examples of subtree perfectness in the literature have been total preorders failing subtree perfectness in exactly these trees, which may have created the impression that only expected utility can be used in such problems, and it may not be so well-known that alternatives such as partial orders can be useful here too.
Problems in which one will observe some data and then make a single decision are not rare. Therefore, it is worth knowing that, in such cases, expected utility does not have a monopoly on subtree perfectness.

In this paper we have only considered a subset of trees of this type: those in Fig.~\ref{fig:subtree:perfecness:simple:trees}\subref{fig:subtree:perfecness:simple:trees:1} and the one in Fig.~\ref{fig:subtree:perfecness:eadm:tree}.
We can consider a wider class of trees of this type: those with any size of partition at the first chance node, with each chance arc leading to a decision node involving any number of gambles.
We could then try to find conditions for subtree perfectness in all trees of this type, and to find how one has to restrict particular choice functions to attain this behavior (for instance, an extension of Example~\ref{ex:gammamaximin:success:chancenodes}).
Some preliminary results can be found in~\cite{2011:huntley:thesis}.

Subtree perfectness has strong links with several other properties of the solutions of sequential decision problems.
The first is the classical method of backward induction: solving decision trees by ``rolling back''.
This algorithm has not often been employed for criteria other than total preorders, but two similar algorithms have been proposed by Huntley and Troffaes~\cite{2012:huntley::backinduct} and Kikuti et al.~\cite{2011:kikuti:sequential}.
In particular the latter reference contains all the information required to find the normal form solution for maximality and E-admissibility via backward induction.
We can ask when backward induction gives the same solution as the normal form solution we have defined in this paper.
As both the above works note, it is possible for backward induction to work but subtree perfectness to fail.
Again, this typically never occurs when using a total preorder, so this result may not be well known.
It is, however, quite easy to show that subtree perfectness implies backward induction~\cite[Corollary~27]{2011:huntley:subtree:perfectness}.

Another related question is the equivalence of our normal form solutions with \emph{extensive form} solutions.
In extensive form solutions, the subject makes decisions only upon reaching a decision node, rather than making all decisions ahead of time.
Such a solution can be seen as a subtree of the initial tree, with some decision arcs removed.
Upon reaching a decision node, the subject then simply picks one of the arcs that remains.
Since subtree perfectness refers to consistency between local and global solutions, one might expect it to be sufficient for the equivalence between the two forms, but consider Example~\ref{ex:eadmissibility:success:special}.
Here, subtree perfectness is satisfied, but there is no extensive form solution available that we could really call equivalent.
This behavior arises because there are many more possible normal form solutions than there are extensive form solutions.
This topic requires further study; introductory material can be found in \cite[\S~7.2]{2011:huntley:subtree:perfectness}.

%
%
\section*{Acknowledgments}

The authors are indebted to Teddy Seidenfeld for suggesting the term `subtree perfectness'. EPSRC supported the first author.

\bibliographystyle{amsplainurl}
\bibliography{imprecisetrees}

\end{document}